\documentclass[12pt]{article}

\usepackage{anysize  }
\usepackage{amssymb  }
\usepackage{amsmath  }
\usepackage{amsthm   }
\usepackage[dvips]{graphicx}
\usepackage{color    }
\usepackage{url      }
\usepackage{paralist }
\usepackage{psfrag   }
\usepackage{wrapfig  }
\usepackage{subfigure}

\newcommand{\conv}{\operatorname{conv}}

\newcommand{\R}{{\mathbb R}}
\newcommand{\Z}{{\mathbb Z}}
\newcommand\ov[1]{\overline{#1}}
\newcommand  {\D   }   {{\mathbb D}}%
\renewcommand{\S   }   {{\mathbb S}}%
\renewcommand{\H   }   {{\mathbb H}}%
\newcommand  {\eps }   {\varepsilon}%

\theoremstyle{plain}
\newtheorem{theorem}{Theorem}[section]
\newtheorem{lemma}[theorem]{Lemma}
\newtheorem{proposition}[theorem]{Proposition}
\newtheorem{corollary}[theorem]{Corollary}

\newtheorem*{lemma*}{Lemma}

\theoremstyle{definition}
\newtheorem{definition}[theorem]{Definition}
\newtheorem{example}[theorem]{Example}

\theoremstyle{remark}
\newtheorem*{remark}{Remark}

\allowdisplaybreaks[4]%

\title{The $E_t$-Construction for Lattices,\\
       Spheres and Polytopes}
\author{Andreas Paffenholz\thanks{%
    Research supported by the Deutsche Forschungsgemeinschaft within
    the European graduate program `Combinatorics, Geometry, and
    Computation' (No. GRK 588/2)} 
\setcounter{footnote}{6}
\qquad G\"unter M. Ziegler\thanks{%
  Partially supported by Deutsche Forschungsgemeinschaft, 
  via the DFG Research Center ``Mathematics in the Key Technologies'' (FZT86), 
  the Research Group ``Algorithms, Structure, Randomness'' (Project ZI 475/3), 
  a Leibniz grant (ZI 475/4), and by the German Israeli Foundation (G.I.F.)}\\
\small Inst.\ Mathematics, MA 6-2\\[-1.6mm]
\small TU Berlin, D-10623 Berlin, Germany\\[-1.6mm]
\small\url{{paffenholz,ziegler}@math.tu-berlin.de}}

\begin{document}
\date{\tiny March 17, 2004} 
\maketitle


\begin{abstract}\noindent%
  We describe and analyze a new construction that produces new
  Eulerian lattices from old ones. It specializes to a construction
  that produces new strongly regular cellular spheres (whose face
  lattices are Eulerian).
  
  The construction does not always specialize to convex polytopes;
  however, in a number of cases where we can realize it, it produces
  interesting classes of polytopes. Thus we produce an infinite family
  of rational $2$-simplicial $2$-simple $4$-polytopes, as requested by
  Eppstein, Kuperberg and Ziegler \cite{Z80}. We also construct for
  each $d\ge3$ an infinite family of $(d-2)$-simplicial $2$-simple
  $d$-polytopes, thus solving a problem of Gr\"unbaum~\cite{Gr1-2}.
\end{abstract}


\section*{Introduction}

Eulerian lattices provide an interesting, entirely combinatorial model
for the combinatorics (face lattices) of convex polytopes. Although
the concepts have been around for a long time (Eulerian posets were
formalized by Stanley \cite{stanley82:_some} in 1982, but the ideas
may be traced back at least to Klee's paper~\cite{Klee3} from 1964),
there are surprisingly little systematic knowledge and development of
Eulerian \emph{lattices}, despite a number of extensive recent studies
of (the flag vectors of) Eulerian \emph{posets};
cf.~Stanley~\cite{St}.  A complete list of linear relations for the
flag vectors of polytopes, spheres and lattices was given by Bayer and
Billera in \cite{BaBi}, where they proved their generalized
Dehn-Sommerville relations.

In this paper we present a combinatorial construction $E_t$ for
Eulerian lattices which specializes to the setting of ``strongly
regular spheres.'' Here the parameter $t$ is an integer between $0$
and $\ell-2$, where $\ell$ is the length of the lattice, and denotes
the dimension (rank plus one) of the elements that will correspond to
coatoms of the new lattice. The $E_t$-construction may be more
intuitive when applied to cellular spheres, but Eulerian lattices
provide a simple axiomatic setting in which we can most easily analyze
the construction and derive its properties. Thus we will start with
lattices (Section~\ref{sec:EulerianLattices}), and in a second step
transfer the results to spheres and give a geometric interpretation of
the construction (Section~\ref{sec:CellularSpheres}).  In some cases
this construction can be performed in the setting of convex polytopes,
and then be used to construct interesting new classes of such
(Section~\ref{sec:Polytopes}):
\begin{compactenum}
\item We provide the first infinite family of \emph{rational}
  $2$-simplicial $2$-simple $4$-polytopes.
\item We present for each $d\ge 3$ an infinite family of
  $(d-2)$-simplicial $2$-simple $d$-polytopes.
\end{compactenum}
The first class solves a problem posed by Eppstein, Kuperberg and
Ziegler (\cite{Z80}, compare also \cite{Z82}). They provided an
infinite family of $2$-simplicial $2$-simple $4$-polytopes, but their
construction is quite rigid and produces nonrational coordinates. In
contrast, our construction will allow great flexibility for the actual
choice of coordinates for a geometric realization of the polytopes in
this family. Thus we can also derive a number of consequences for the
$f$-vector classification problem for $4$-dimensional polytopes
(Section~\ref{sec:fVectors}). The second class of polytopes solves a
problem of Gr\"unbaum \cite[Problem 9.7.7(iii)]{Gr1-2} and thus
substantiates claims by Perles and Shephard \cite[p.170]{Gr1-2} 
and by Kalai \cite[19.5.19]{Kal97}.

Our construction contains the method of Eppstein, Kuperberg and
Ziegler as a special case: If we choose $t=1$ and apply $E_1$ to the
face lattice of a simplicial $4$-polytope $P$ whose edges are tangent
to the unit sphere, then we obtain the face lattice of a
$2$-simplicial $2$-simple polytope which they call $E(P)$. Thus the
second class of examples displayed above is obtained by a
generalization of the approach from \cite{Z80} to higher dimensions.

In contrast, the construction of the first class relies on the
analysis of the dual construction to $E_1$, which we describe as
``vertex truncation'' in Section~\ref{subsec:vt}.  The special case
where $P$ is a regular polytope appears in the literature in the
context of the construction of regular and ``uniform'' polytopes: Here
our vertex truncation can be seen as an instance of ``Wythoff's
construction'' as described by Coxeter (\cite[p.~147]{Cox} and
\cite{Cox2}), while a special case of the dual $E_t$-construction
appears in G\'evay's work \cite{Gevay,Gevay2}.  
\smallskip

\emph{Acknowledgement. } We would like to thank Carsten Lange, Julian
Pfeif\-le and Arnold Wa\ss mer for several discussions and valuable
hints. 


\section{Eulerian lattices}\label{sec:EulerianLattices}

We refer to Stanley \cite{Sta5} for general terminology and background
about posets, lattices, and their enumerative combinatorics.
All the posets we consider are finite. A \emph{bounded} poset is
required to have a unique maximal element $\hat1$ and a unique minimal
element $\hat0$.  A poset $L$ is a \emph{lattice} if any two elements
$y,y'\in L$ have a unique minimal upper bound $y\vee y'$ and a unique
maximal lower bound $y\wedge y'$, known as the \emph{join} and the
\emph{meet} of $y$ and~$y'$, respectively. By $\bigvee A$ and
$\bigwedge A$ we denote the join resp.\ the meet of a subset
$A\subseteq L$. A lattice is \emph{complemented} if for every element
$x$ there is an element $x'$ with $x\wedge x'=\hat0$ and $x\vee
x'=\hat1$.

A bounded poset is \emph{graded} if all its maximal chains have the
same length. Every graded poset comes with a \emph{rank function}
$r:L\rightarrow\Z$, normalized such that $r(\hat0)=0$.  The
\emph{length} of a graded poset is given by~$r(\hat1)$.  By $L_i$ we
denote the set of all elements of rank~$i+1$ in~$L$, for $-1\le i \le
r(\hat 1)-1$. We also talk about the \emph{dimension}
$\dim(x):=r(x)-1$ of an element $x\in L$: This is motivated by the
important situation when $L$ is the face lattice of a polytope, and
$x\in L_i$ corresponds to a face of dimension~$i$. In this case we
write $d:=r(\hat 1)-1$ for the dimension of the polytope.

By $f_i(L)=|L_i|$ we denote the number of elements of rank $i+1$ (for
$-1\le i\le d$), while $f_{ij}(L)$, for $-1\le i\le j\le d$, denotes
the number of pairs of elements $x\in L_i$ and $y\in L_j$ with $x\le
y$. If $L$ is a bounded poset, we denote by
$\ov{L}:=L{\setminus}\{\hat0,\hat1\}$ its \emph{proper part}.  The
minimal elements of $\ov{L}$, i.\,e.\ the elements of $L$ of rank~$1$,
are the \emph{atoms} of $L$; similarly the maximal elements of
$\ov{L}$, i.\,e.\ the elements of $L$ of rank~$d$, are the
\emph{coatoms} of $L$.  By $L^{op}$ we denote the \emph{opposite} of a
poset $L$, that is, the same set, but with reversed order relation.
\begin{definition}
  A graded poset $L$ of length $d+1$ is \emph{Eulerian} if every
  interval $[x,y]$ with $x<y$ has the same number of elements of odd
  and of even rank.
\end{definition}
This definition is equivalent to the ``usual'' definition
\cite[p.~122]{Sta5}, which requires that $\mu(x,y)=(-1)^{r(y)-r(x)}$
for all $x\le y$, where $\mu$ is the M\"obius function
\cite[Ex.~3.69(a)]{Sta5}.
\begin{definition}[The $\mathbf E_t$-construction for Eulerian posets]
  Let $L$ be a graded poset of length $d+1$, and let $t$ be any
  integer in the range $0\le t\le d-1$.
  
  The \emph{$E_t$-construction} applied to $L$ yields the poset
  $E_t(L)$ that consists of the following subsets of~$L$, ordered by
  reversed inclusion:
  \begin{compactitem}[~$\bullet$]
  \item the empty set,
  \item the one element sets $\{y\}$ for $y\in L_t$, and
  \item the intervals $[x,z]\subseteq L$ such that some $y\in L_t$
    satisfies $x<y<z$.
  \end{compactitem}
\end{definition}
\begin{figure}[htb]
  \begin{center}
    \psfrag{d+1}{$\displaystyle d+1$} \psfrag{t+1}{$\displaystyle
      t+1$} \psfrag{0}{$\displaystyle 0$} \psfrag{1}{$\displaystyle
      1$} \psfrag{x}{$\displaystyle x$} \psfrag{z}{$\displaystyle z$}
    \psfrag{L}{$L$} \psfrag{E_t(L)}{$E_t(L)$}
    \psfrag{[x,z]}{$\displaystyle [x,z]$}
    \psfrag{emptyset}{$\displaystyle \emptyset$}
    \psfrag{r(x)+d+1-r(z)}{$\displaystyle r(x)+d+1-r(z)$}
    \psfrag{[x,1]}{$\displaystyle [\text{atom},\hat 1]$}
    \psfrag{[0,z]}{$\displaystyle [\hat 0,\text{coatom}]$}
    \psfrag{atom}{$\displaystyle\text{atom}$}
    \psfrag{coatom}{$\displaystyle\text{coatom}$}
    \psfrag{{y}}{$\displaystyle \{y\}$} \psfrag{d}{$\displaystyle d$}
    \includegraphics[height=4.5cm]{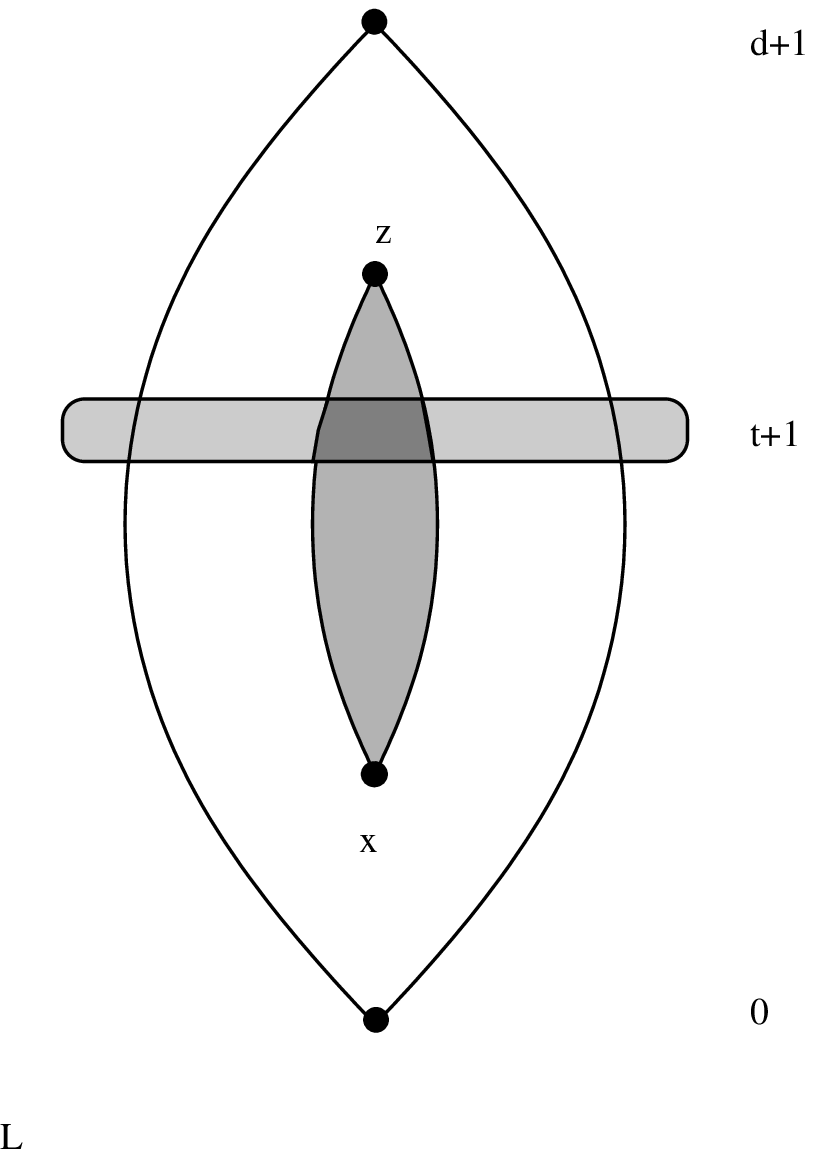} \hspace{2cm}
    \includegraphics[height=4.5cm]{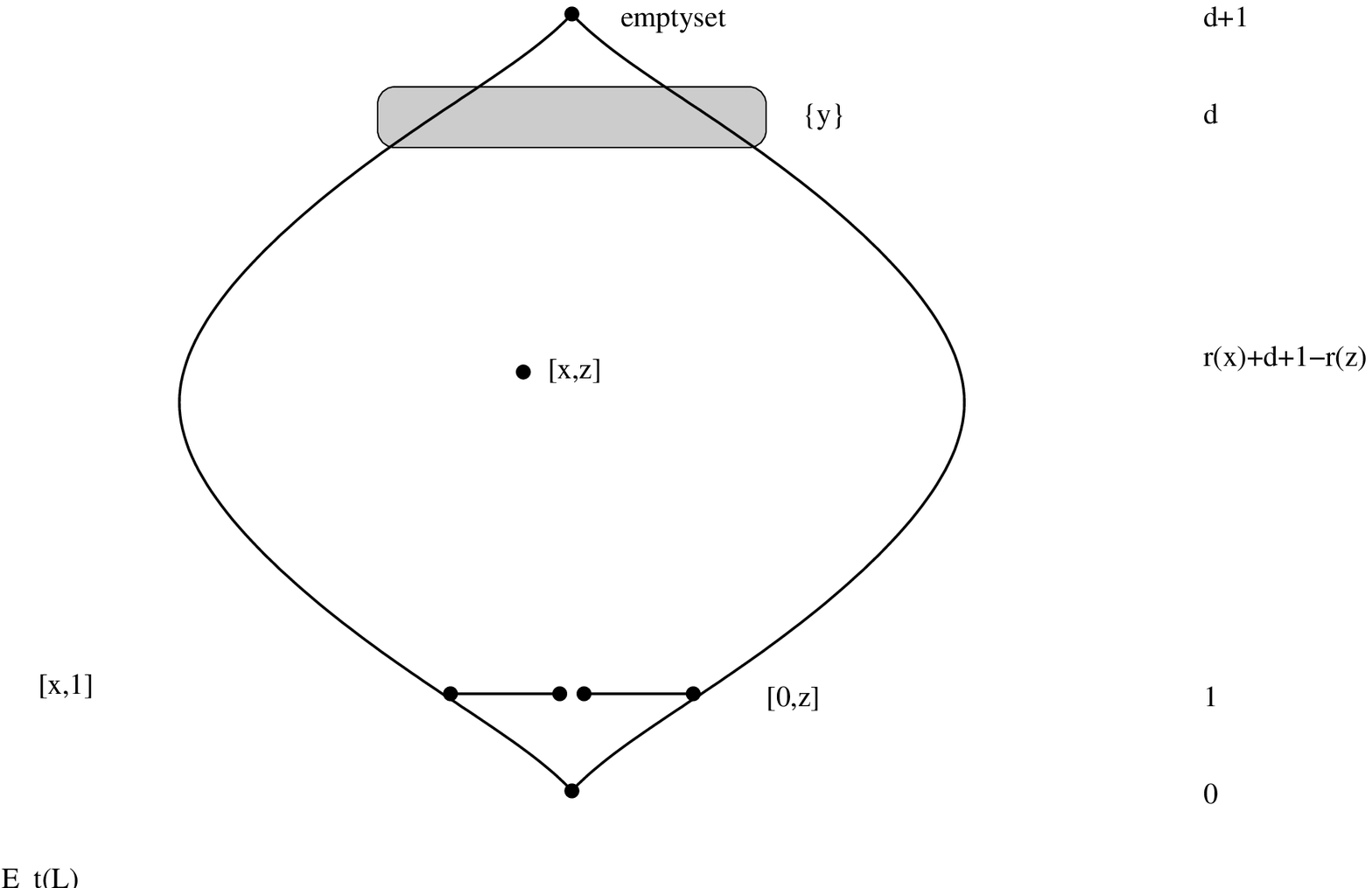}
    \caption{Combinatorial construction of $E_t(L)$}
    \label{facelattice}
  \end{center}
\end{figure}
The rank function on $E_t(L)$ is given by
\begin{eqnarray*}
\rho(b)=
  \begin{cases}
    r(x)+d+1-r(z) \qquad \textrm{ for } & b=[x,z],\quad r(x)<t+1<r(z),\\
    d & b=\{y\},\quad y\in L_t\\
    d+1 & b=\emptyset.
  \end{cases}
\end{eqnarray*}
Thus $E_t(L)$ is again a graded poset of length $d+1$. Its coatoms are
the one-element sets $\{y\}$, $y\in L_t$, its atoms are the intervals
$[x,\hat1]$ (for $1\le t\le d-1$) and $[\hat0,z]$ (for $0\le t\le
d-2$), where the $x$ and $z$ are the atoms resp.\ coatoms of~$L$. The
$f$-vector is
\begin{equation}\label{eq:fvector}
f_k(E_t(L))\ \ =\ \ 
  \begin{cases}
    \sum_{i,j} f_{ij}(L) & \textrm{ for }-1\le k<d-1,\\
    f_t(L)               & \textrm{ for }k=d-1,\\
    1 & \textrm{ for }k=d,
  \end{cases}
\end{equation}
where the sum in the above formula is over all pairs $i,j$ in the
range $-1\le i<t<j\le d$ with $j-i=d-k$.

In the following we shall mostly be interested in \emph{Eulerian
  lattices}, that is, in Eulerian posets which at the same time
satisfy the lattice property.  If $L$ is a lattice, then the set
$E_t(L)$ can more compactly be described as
\[\textstyle
E_t(L)\ \ :=\ \ \left\{[\bigwedge A,\bigvee A]:A\subseteq L_t\right\},
\]
again ordered by reversed inclusion.  In this case $E_t(L)$ is again a
lattice: The join and meet operations in $E_t(L)$ are given by
\[
[x,z]\wedge[x',z']\ =\ [x\wedge x',z\vee z'],\qquad [x,z]\vee [x',z']\ 
=\ [x\vee x',z\wedge z'].
\]
This follows from the fact that every interval $[x,z]$ is itself an
Eulerian lattice, hence it has non-zero M\"obius number, hence it is
complemented, and the join of its atoms is the maximal element and the
meet of the coatoms is the minimal element (see
\cite[Cor.~3.9.5]{Sta5}).

\begin{remark}
  Even if $L$ is an Eulerian poset that is not a lattice, $E_t(L)$ may
  still be a lattice.
\end{remark}

\begin{example}
  Let $L$ be a graded poset of length $d+1$.  For $t=0$ we obtain
  $E_0(L)\cong L^{op}$.  For $t=d-1$ we get $E_{d-1}(L)\cong L$. Thus
  the ``boundary cases'' of $t=0$ and $t=d-1$ are not interesting, and
  hence they will be excluded from some of the following
  discussions without loss of interesting effects.  We also note that 
  \[
  E_t(L)\ =\ E_{d-1-t}(L^{op}),
  \]
  so we can derive the same posets/lattices from $L$ and from
  $L^{op}$.
\end{example}

\begin{theorem}\label{thm:eulerian}
  Let $L$ be an Eulerian poset of length $d+1$ and $0\le t\le d-1$.
  Then $E_t(L)$ is an Eulerian poset of the same length~$d+1$.
\end{theorem}

\begin{proof}
  This is true for $t=0$ and for $t=d-1$, which includes all possible
  $t$-values for $d\le2$.  Thus we may use induction on the length
  of~$L$.
  
  First we show that all proper intervals in $E_t(L)$ are Eulerian.
  Indeed, for any element $[x,z]\in E_t(L)$, $(x,z)\ne(\hat0,\hat1)$,
  the upper interval $\big[[x,z],\hat1\big]$ of $E_t(L)$ is isomorphic
  to $E_{t'}([x,z])$ for $t'=t-r(x)$; hence all proper upper intervals
  in $E_t(L)$ are produced by the $E_t$-construction from Eulerian
  posets of smaller length, so they are Eulerian by induction.
  
  Similarly, if $[x,z]$ is an element of rank at most $d-1$ in
  $E_t(L)$, that is, an interval of $L$ with $x<y<z$ for some $y\in
  L_t$, then the lower interval $\big[\hat0,[x,z]\big]$ of $E_t(L)$ is
  isomorphic to $[\hat0,x]\times[z,\hat1]^{op}$, which is Eulerian
  since $L$ is Eulerian. If $\{y\}$ is a coatom of $E_t(L)$, for $y\in
  L_t$, then the lower interval $\big[\hat0,\{y\}\big]$ is isomorphic
  to $(L_{<y})\times (L_{>y})^{op}\ \uplus\ \hat1$. Thus it is a
  ``reduced product'' and this operation is known
  \cite[Ex.~3.69(d)]{Sta10} and easily checked to preserve the
  Eulerian property; compare Walker \cite[Sect.~6]{Walker-canon}.
  
  Finally, we have to see that $E_t(L)$ itself has the same number of
  odd and even rank elements. For this we may use the $f$-vector of
  $E_t(L)$, as computed above.  Every interval $[\hat 0,z]$ is
  Eulerian. Thus for $0\le j\le d-1$ and all $z\in L_j$ we have
  $\sum_{i=-1}^j(-1)^if_i([\hat 0,z])=0$, which by summing over all
  $z\in L_j$ yields $\sum_{i=-1}^j(-1)^if_{ij}=0$, that is,
  \begin{align}
    \sum_{i=-1}^{t-1}(-1)^if_{ij}&=\;-\sum_{i=t}^j(-1)^if_{ij},\label{1}
    \intertext{for $j\ge t\ge0$.  This is one of the generalized
      Dehn-Sommerville equations \cite{BaBi}. A similar argument for
      upper intervals shows that}\notag\\[-10mm]
    \sum_{j=i}^d(-1)^jf_{ij}&=\;(-1)^d\delta_{id}\label{2}.
  \end{align}
  for $i\le d$.  With these two equations, we can compute
  \begin{alignat*}{3}
    \sum_{k=-1}^d(-1)^{d-k}&f_k(E_t(L))\ \ =\\[-3mm]
     =&\hspace{12pt}1-f_t+\sum_{i=-1}^{t-1}\sum_{j=t+1}^d(-1)^{j-i}f_{ij}
    &=&\hspace{12pt}1-f_t+\sum_{j=t+1}^d(-1)^j\sum_{i=-1}^{t-1}(-1)^if_{ij}\\
    {\buildrel(\ref{1})\over=}&\hspace{12pt}
    1-f_t-\sum_{j=t+1}^d(-1)^j\sum_{i=t}^{j}(-1)^if_{ij}\quad
    &=&\hspace{12pt}1-\sum_{j=t}^d(-1)^j\sum_{i=t}^{j}(-1)^if_{ij}\\
    =&\hspace{12pt} 1-\sum_{i=t}^d(-1)^i\sum_{j=i}^{t}(-1)^jf_{ij}
    &{\buildrel(\ref{2})\over=}&\hspace{12pt}
    1-\sum_{i=t}^d(-1)^i(-1)^d\delta_{id}\qquad =\quad 0.
  \end{alignat*}
  
  Alternatively, one may argue from Theorem~\ref{thm:homeo} in the
  next section: Since the order complexes of $L$ and of $E_t(L)$ are
  homeomorphic, they must have the same Euler characteristic, which is
  the M\"obius function of $L$ resp.\ $E_t(L)$, which is what we need
  for $(x,z)=(\hat 0,\hat 1)$.
\end{proof}

\begin{definition} 
  Let $L$ be an Eulerian lattice of length $d+1$, and let $0\le k,h\le
  d-1$.
  \begin{compactitem}[~$\bullet$]
  \item $L$ is \emph{boolean} if it is isomorphic to
      the face lattice of a $d$-simplex.
    \item $L$ is \emph{$k$-simplicial} if all intervals $[\hat 0,z]$
      with $r(z)=k+1$ are boolean.
    \item It is \emph{$h$-simple} if all intervals $[x, \hat 1]$ with
      $r(x)=d-h$ are boolean.
    \item $L$ is {\em simple} resp.\ {\em simplicial} if it is
      $(d-1)$-simple resp.\ $(d-1)$-simplicial.
  \end{compactitem}
  A \emph{$(k,h)$-lattice} $L$ is a $k$-simplicial and $h$-simple
  Eulerian lattice.
\end{definition}

Every $k$-simplicial Eulerian lattice is also $(k-1)$-simplicial for
$k>0$. All Eulerian lattices are $1$-simplicial.  The property
``$k$-simplicial'' is dual to ``$k$-simple''; thus every Eulerian
lattice is also $1$-simple.  $L$ is a $(k,h)$-lattice if and only if
$L^{op}$ is an $(h,k)$-lattice.

\begin{remark}
  Every $(k,h)$-lattice of length $d+1$ with $k+h>d$ is boolean.  (In
  particular, any $k$-simplicial $h$-simple $d$-polytope with $k+h>d$
  is a simplex.) The proof, also in the generality of Eulerian
  lattices, is a straightforward extension of the argument for simple
  and simplicial polytopes (see \cite[p.~65]{Gr1-2} and
  \cite[p.~23]{Ziegler95}). 
\end{remark}

Any Eulerian lattice is graded and complemented, so by induction on
the length it follows that every Eulerian lattice $L$ of length $\ell$
has at least $\binom{\ell}i$ elements of rank~$i$. Furthermore, if
equality holds for some $i$, $0<i<\ell$, then $L$ is boolean.  In
particular, any Eulerian lattice of length $\ell$ has at least $\ell$
atoms, with equality only if $L$ is boolean, and similarly for
coatoms.  We rely on this criterion in the proof of the following
characterization of Eulerian lattices $L$ for which $E_t(L)$ is
$k$-simplicial and $h$-simple.

\begin{theorem}
  Let $L$ be an Eulerian lattice of length $d+1$ and let
  $1\le t\le d-2$.
  \begin{compactitem}
  \item[\rm(1)] For $0\le k\le d-2$, the lattice $E_t(L)$ is
    {$k$-simplicial} if and only if $L$ is $s$-simplicial and
    $r$-simple for $s=\min\{k,t-1\}$ and $r=\min\{k,d-t-2\}$.
  \item[\rm(1$'$)] The lattice $E_t(L)$ is never $(d-1)$-simplicial.
  \item[\rm(2)] The lattice $E_t(L)$ is {$2$-simple} if and only if
    every interval $[\overline{x},\overline{z}]\subseteq L$
    with $r(\overline{x})=t-1$ and $r(\overline{z})=t+3$ is boolean.\\
    In particular, this is the case if $L$ is $(t+2)$-simplicial or
    $(d-t+1)$-simple. 
  \item[\rm(2$'$)] The lattice $E_t(L)$ is never $3$-simple.
  \end{compactitem}
\end{theorem}

\begin{proof}
  (1): The elements of rank at most $k+1\le d-1$ in $E_t(L)$ are the
  intervals of the form $a=[x,z]\subseteq L$ with $r(x)\le t$,
  $r(z)\ge t+2$ satisfying $r(x)+(d+1)-r(z)=\rho([x,z])\le k+1$.  An
  element $x\in L$ appears as the lower end of such an interval
  $[x,z]$ if and only if $0\le r(x)\le\min\{k+1,t\}$;
  similarly this concerns all elements $z\in L$ of corank
  $0\le d+1-r(z)\le\min\{k+1,d-t-1\}$.
  
  The atoms of $E_t(L)$ below $a=[x,z]$ are given by both the atoms of
  $L$ below $x$, whose number is at least $r(x)$, and the coatoms of
  $L$ above $z$, whose number is at least $d+1-r(z)$.  So the interval
  has at least $r(x)+d+1-r(z)=k+1$ atoms, with equality if and only if
  the intervals $[\hat0,x]$ and $[z,\hat1]$ of $L$ are both boolean.
  
  Thus \emph{all} lower intervals $[\hat0,a]$ of rank $k+1$ in
  $E_t(L)$ are boolean if and only if \emph{all} intervals $[\hat0,x]$
  and $[z,\hat1]$ are boolean for $r(x)\le \min\{k+1,t\}$ resp.\ 
  $d+1-r(z)\le\min\{k+1,d-t-1\}$.
   
  (1$'$): An analysis as for (1) shows that for any element $\{y\}$ of
  rank $d$ in $E_t(L)$, that is for $y\in L_t$, there are at least
  $t+1$ atoms in $L$ below $y$ and at least $d-t$ coatoms in $L$
  above~$y$; thus there are at least $(t+1)+(d-t)=d+1$ atoms below
  $\{y\}$ in $E_t(L)$: too many.
   
  (2): $E_t(L)$ is $2$-simple if all intervals $[b,\hat1]\subset
  E_t(L)$ with $b=[x,z]\subset L$, $\rho(b)=d-2$, are boolean, that
  is, they have $3$ atoms or coatoms. This is the case if and only if
  every interval $b=[x,z]\subset L$, $r(x)<t+1<r(z)$, of length $3$
  contains precisely three elements of rank $t+1$. This is equivalent
  to the condition that every length~$4$ interval $[\ov x, \ov z]$
  with $r(\ov x)=t-1$ and $r(\ov z)=t+3$ is boolean.  In terms of the
  usual flag vector notation, this can numerically be expressed as
  $f_{t-2,t,t+2}(L)=6f_{t-2,t+2}(L)$.
  
  (2$'$): Similarly, for $E_t(L)$ to be $3$-simple we would need that
  every interval $[x,z]$ in $L$ of length $4$ with $r(x)<t+1<r(z)$ and
  $r(z)=r(x)+4$ contains exactly $4$ elements of rank $t+1$. This is
  impossible for the case where $r(x)=t-1$ and $r(z)=t+3$, where the
  interval $[x,z]$ has at least $6$ elements in its ``middle level''
  (that is, of rank~$t+1$) by the Eulerian condition.
\end{proof}

\begin{corollary}\label{cor:2s2s:lattices}
  For $d\ge3$ and any simplicial Eulerian lattice $L$ of length $d+1$
  the lattice $E_{d-3}(L)$ is $(d-2)$-simplicial and $2$-simple.\qed
\end{corollary}


\section{CW spheres}\label{sec:CellularSpheres}

The \emph{order complex} of a bounded poset $L$ is the abstract
simplicial complex of chains in its proper part ~$\ov{L}$, denoted by
$\Delta(\ov{L})$. By $\|\ov{L}\|$ we denote the geometric realization
of the order complex of the poset $\ov L$.

A CW complex is \emph{regular} if all its closed cells are embedded,
that is, if the attaching maps of the individual cells make no
identifications on the boundary.  In the following, all \emph{cell
  complexes} are finite regular CW complexes. A cell complex is
\emph{strongly regular} if the following ``intersection property''
\cite[Problem 4.47, p.~223]{Z10} holds: The intersection of any two
cells is a single closed cell (which may be empty). This holds if and
only if not only the cells, but also all the stars of vertices are
embedded.  For example, the boundary complex of any convex
$d$-polytope is a strongly regular CW $(d-1)$-sphere.

The main property of PL spheres and manifolds (see \cite{Hudson} for
the definitions) we need is the following: They come with regular cell
decompositions such that not only the boundary of each cell is a
sphere (of the appropriate dimension), but also the links of all faces
are genuine spheres (not only homology spheres, as in general
triangulated or cellular manifolds \cite[Prop.~4.7.21]{Z10}). 
Equivalently, in the face poset of the cell decomposition, augmented
by a maximal element $\hat1$, not only the order complexes of the
lower intervals $[\hat0,x]$ with $x<\hat1$ are spheres, but the same
is true for \emph{all} intervals $[x,y]$, with the only possible
exception of $[\hat0,\hat1]$, whose order complex is homeomorphic to
the manifold in question.

Finally, a bounded, graded poset $P$ is the face poset of a regular CW
sphere if and only if the order complex of every interval $[\hat0,x]$,
$x\in P$, is homeomorphic to a sphere of dimension $r(x)-2$; see
\cite[Prop.~4.7.23]{Z10}.

The following is similar to (simpler) results and proofs for interval
posets in Walker~\cite{Walker-canon}.

\begin{theorem}\label{thm:homeo}
  Let $L$ be a graded poset of length $d+1$. Then $L$ and $E_t(L)$ are
  PL-homeomorphic:
  \[
  \|\ov{L}\| \ \cong\ \|\ov{E_t(L)}\|.
  \]
\end{theorem}

\begin{figure}[h]\label{fig:subdiv}
\begin{center}
  \input{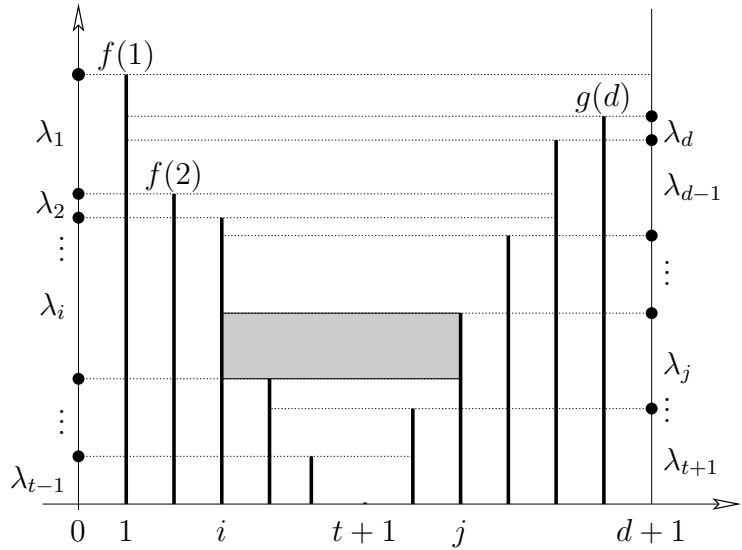}
\end{center}
\caption{Sketch for the proof of Theorem \ref{thm:homeo}.
  The height of the shaded rectangle indicates the size of the
  coefficient $\alpha_{i,j}$.}
\end{figure}

\begin{proof}
  We verify that $\Delta\bigl(\ov{E_t(L)}\bigr)$ is a subdivision
  of~$\Delta\bigl(\ov{L}\bigr)$, and give explicit formulas for the
  subdivision map and its inverse. (Compare this to Walker
  \cite[Sects.~4,5]{Walker-canon}.)
  
  A canonical map
  $\pi:\|\Delta\bigl(\ov{E_t(L)}\bigr)\|
                \longrightarrow\|\Delta\bigl(\ov{L}\bigr)\|$,
  linear on the simplices, is given by its values on the vertices,
  \[\begin{array}{rcl@{\quad}l}
    \{y\}    &\longmapsto&  y               & \textrm{for }y\in L_t, \\
  {}[x,z]    &\longmapsto&\frac12x+\frac12z & \textrm{for }
                                    \hat0<x<y<z<\hat1, \ y\in L_t, \\
  {}[x,\hat1]&\longmapsto&  x               & \textrm{for }
                                    \hat0<x<y<\hat1, \ y\in L_t, \\
  {}[\hat0,z]&\longmapsto&  z               & \textrm{for }
                                    \hat0<y<z<\hat1, \ y\in L_t.
  \end{array}\]
  This map is well-defined and continuous. Its inverse, a subdivision
  map, may be described as follows: Any point of
  $\Delta\bigl(\ov{L}\bigr)$ is an affine combination of elements on a
  chain in~$\ov{L}$, so it may be written as
  \[
  \mathbf{x}: \lambda_1 x_1 \ <\ \ \cdots\ \ <\ \lambda_t x_t \ <\ 
  \lambda_{t+1} y_{t+1} \ <\ \lambda_{t+2} z_{t+2}\ <\ \ \cdots\ \ <\ 
  \lambda_d z_d.
  \]
  with $\lambda_i\ge0$ and $\sum_i\lambda_i=1$.  We set $x_0:=\hat0$
  and $z_{d+1}:=\hat1$, with $\lambda_0:=1$ and $\lambda_{d+1}:=1$.
  Now the above point is mapped by $\pi^{-1}$ to
  \[
  \pi^{-1}(\mathbf{x})\ \ =\ \ \lambda_{t+1}~\{y_{t+1}\}\ +\ 
  \sum_{1\le i< t,\, t+1<j\le d}2\alpha_{i,j}~[x_i,z_j] \ +\ 
  \sum_{i=0,\, t+1<j\le d \textrm{ or}\atop 1\le i\le
    t,\,j=d+1}\alpha_{i,j}~[x_i,z_j],
  \]
  where the coefficients $\alpha_{i,j}$ are given by
  \begin{align*}
    \alpha_{i,j} & =
    \begin{cases}
      \min\{f(i), g(j)\} \ -\ \max\{f(i+1),g(j-1)\}
      &\quad \textrm{if this is }\ge0,\\
      0 &\quad \textrm{otherwise},
    \end{cases}
  \end{align*}
  with
\[  \begin{array}{lll}
    f(i) & := \ \lambda_i+\lambda_{i+1}+\ldots+\lambda_t 
            &\textrm{ for } 0\le i\le t+1\textrm{, and }\\
    g(j) & := \ \lambda_{t+2}+\ldots+\lambda_{j-1}+\lambda_j 
            &\textrm{ for } t+1\le j\le d+1
  \end{array}
\]
(with $f(t+1)=g(t+1)=0$). We refer to Figure \ref{fig:subdiv} for
illustration.
\end{proof}

\begin{theorem}\label{thm:2.2}
  If $L$ is the face poset of a regular CW PL sphere or PL manifold,
  then so is $E_t(L)$.
\end{theorem}

\begin{proof}
  By Theorem~\ref{thm:eulerian} and its proof, using the PL property,
  we get the cell complex.  By Theorem~\ref{thm:homeo}, this cell
  complex is homeomorphic to $\|\ov L\|$.
\end{proof}

The following is an easy consequence of Theorem \ref{thm:2.2} and
Corollary \ref{cor:2s2s:lattices}.

\begin{corollary}\label{cor:2s2s}
  For $d\ge3$ and for any strongly regular simplicial PL-sphere $S$ of
  dimension $d-1$, the PL-sphere $E_{d-3}(S)$ is $(d-2)$-simplicial
  and $2$-simple.\qed
\end{corollary}


\section{Polytopes}\label{sec:Polytopes}

We refer to \cite{Gr1-2} and \cite{Ziegler95} for background on
polytopes.  The boundary of any polytope $P$ naturally carries the
structure of a strongly regular CW PL-sphere. Thus we can apply the
$E_t$-construction to $\partial P$ and get a new PL-sphere $E_t(P)$.
As mentioned earlier, it is not at all clear that $E_t(P)$ can be
realized as a polytope, for any given convex polytope $P$ and for
given $t$.  However, in the main part of this section we present and
analyze two techniques that do yield infinite families of interesting
$E_t$-polytopes:
\begin{compactitem}[ $\bullet$ ]
\item The first construction is surprisingly simple; it produces the
  first infinite families of rational $2$-simplicial $2$-simple
  $4$-polytopes.
\item The second construction is an extension of the $E$-construction
  of Eppstein, Kuperberg, and Ziegler \cite{Z80} to higher dimensions;
  it produces for each $d\ge3$ infinitely many $(d-2)$-simplicial
  $2$-simple $d$-polytopes; as far as we can see this is new for
  $d>4$, and confirms assertions of G\"unbaum
  \cite[pp.~169-170]{Gr1-2} and Kalai \cite[19.5.19]{Kal97}.
\end{compactitem}
Furthermore, we then survey the previously known examples of
$(d-2)$-simplicial $2$-simple polytopes, and demonstrate that
virtually all of them can be realized as instances of the
$E_t$-construction or of its dual.


\subsection{Realizations via vertex truncation}\label{subsec:vt}

The following simple construction --- if it can be realized ---
produces a polytope that is dual to $E_1(P)$.  A very special case of
it turns up already in Gosset's 1897 construction of the $24$-cell. A
slightly more general version appears in Coxeter's book in the
construction of the three special $4$-dimensional regular polytopes
\cite[pp.~145-164]{Cox}: He considers ``truncation'' of polytopes by
hyperplanes defined by all vertex vectors of a centered regular
$d$-polytope and thinks of this as a continuous process having $d-1$
interesting stages in which the cutting hyperplanes intersect in the
centroids of the $k$-faces, for $1\le k\le d-1$.

\begin{definition}[Vertex Truncation]
  Let $P$ be any $d$-polytope, $d\ge3$.  A vertex $v$ of $P$ is
  \emph{truncated} if we intersect $P$ with a closed half-space that
  does not contain $v$, but contains all other vertices of $P$ in its
  interior.
\end{definition}

Vertex truncation (that is, ``truncating a single vertex'') results
in a new polytope with one additional facet, and with new vertices
$u_e$ corresponding to the edges of $P$ that contain the vertex $v$
of~$P$ that has been truncated.

\begin{definition}[Truncatable Polytopes]
  A $d$-polytope $P$ is \emph{truncatable} if all its vertices can be
  truncated simultaneously in such a way that one single (relative
  interior) point remains from each edge of $P$.  The resulting
  polytope from such a construction is denoted $D_1(P)$.
\end{definition}

Whether this construction can be performed depends on the realization
of $P$, and on a careful choice of the hyperplanes/halfspaces that
truncate the vertices of $P$ --- see the case when $P$ is a possibly
non-regular octahedron. If it can be performed, the resulting polytope
is usually not uniquely determined --- see the case when $P$ is a
simplex.  However, the following proposition establishes that if the
construction can be performed, then the combinatorial type of $D_1(P)$
is uniquely determined, and that it is of interest for our
investigations.

\begin{proposition}
  Let $P$ be any $d$-polytope, $d\ge3$. If $P$ is truncatable, then
  the dual $D_1(P)^*$ of the resulting polytope realizes~$E_1(P)$.
\end{proposition}

\begin{proof}
  The polytope $D_1(P)$ has two types of facets: First there are the
  facets $F'$ that are obtained by vertex truncation from the facets
  $F$ of $P$. (Here we use that $d\ge3$.)  Secondly, there are the
  ``new'' facets $F_v$ that result from truncating the vertices of
  $P$.  The intersection of any two new facets $F_v$ and $F_w$ is
  empty if $v$ and $w$ are not adjacent in~$P$. Otherwise the
  intersection $F_v\cap F_w$ consists of one single point, the new
  vertex $u_e$ given by the edge $e=(v,w)$.  The vertices of $D_1(P)$
  are given by $u_e$, in bijection with the edges $e$ of~$P$.
 
  Thus to check that $D_1(P)$ has the right combinatorics, we only
  check the correct vertex--facet incidences: $u_e$ lies on $F'$ if
  and only if $e$ is an edge of~$F$.  And $u_e$ lies on $F_v$ if and
  only if $e$ is adjacent to~$v$.
  This amounts to the reversed atom--coatom incidences of $E_1(P)$.
\end{proof}

For the goal of our considerations --- producing $(d-2)$-simplicial
$2$-simple $d$-polytopes --- this construction is only interesting in
dimensions $d=3$ and $d=4$, where we can apply it to simplicial
polytopes, trying to get $2$-simple and (for $d=4$) $2$-simplicial
polytopes. In all other dimensions we would have to apply it to a
$3$-simplicial $(d-3)$-simple $d$-polytope to obtain a
$(d-2)$-simplicial and $2$-simple one, and for large $d$ we know
only two such examples in every dimension (the simplex, and the
half-cube $N^d$; see Section~\ref{subsec:KnownExamples}).  So from now
on we are mainly interested in classes of $3$- and $4$-dimensional
polytopes.

A first obvious class of truncatable polytopes is the following.

\begin{proposition}\label{prop:edge-tangent}
  If a $d$-polytope $P$, $d\ge3$, is realized such that all edges are
  tangent to a $(d-1)$-dimensional sphere $S$, then it is truncatable.
\end{proposition}

\begin{proof}
  Consider any vertex $v$ of the polytope $P$ and consider the tangent
  cone $C$ to $S$ with apex $v$. Then the intersection of $C$ with $S$
  is a $(d-2)$-dimensional sphere $S'$. As all edges adjacent to $v$
  are tangent to $S$, they are contained in $C$ and the point of
  tangency must be contained in $S'$. Thus we can cut with the
  hyperplane defined by $S'$.
\end{proof}

By the Koebe--Andreev--Thurston circle packing theorem
\cite[Thm.~4.12]{Ziegler95}, every $3$-polytope has a realization that
is edge-tangent and hence truncatable. In dimension $4$, edge-tangent
polytopes $P$ were constructed and $E_1(P)=D_1(P)^*$ obtained by
Eppstein, Kuperberg, and Ziegler \cite{Z80}, so our construction
subsumes all their examples.  However, we now show that an
``edge-tangency'' is not necessary for the realization of $D_1(P)$ or
$E_1(P)$.

A \emph{stacked} polytope \cite{Kal97} is obtained from a simplex by
repeatedly adding a new vertex ``beyond'' a facet, that is, by glueing
a simplex onto a facet.
 
By construction, any stacked $d$-polytope is simplicial and the vertex
added last is adjacent to precisely $d$ edges.

\begin{theorem}\label{thm:stacktrunc}
  Every stacked $d$-polytope $(d\ge3)$ has a geometric realization
  that is truncatable.
\end{theorem}

\begin{proof}
  Let $P^d_n$ denote a stacked $d$-polytope ($d\ge3$) with $n+d+1$
  vertices, i.\,e.\ a $d$-simplex that has been stacked $n$ times,
  $n\ge0$. We do not care about the actual choice of facets above
  which we have placed a vertex, so that $P^d_n$ might denote many
  combinatorially distinct polytopes. Obviously $P^d_0$ is
  truncatable. Now let $P^d_n$ be any truncatable stacked polytope and
  choose a facet $F$. As $P^d_n$ is simplicial there are precisely $d$
  cutting hyperplanes $H_i$, $1\le i\le d$, that intersect this
  facet~$F$. Place the new vertex $v$ beyond $F'=F\cap D_1(P_n^d)$ and
  such that it is on the same side of $H_i$ as $D_1(P^d_n)$ for all
  $1\le i\le d$. That is, choose $v$ beyond $F'$ and beneath all other
  facets of~$D_1(P^d_n)$.
  When we form the convex hull of $P^d_n$ and $v$, the $d$ new edges
  adjacent to $v$ will intersect all $d$ hyperplanes $H_i$. We can
  choose the hyperplane through these $d$ points as the cutting
  hyperplane for $v$. This gives a realization of $D_1(P^d_{n+1})$.
\end{proof}

\begin{figure}[htb]
  \hbox{ \hspace{-12mm} \includegraphics[height=4cm,bb=66 120 488
    371,clip]{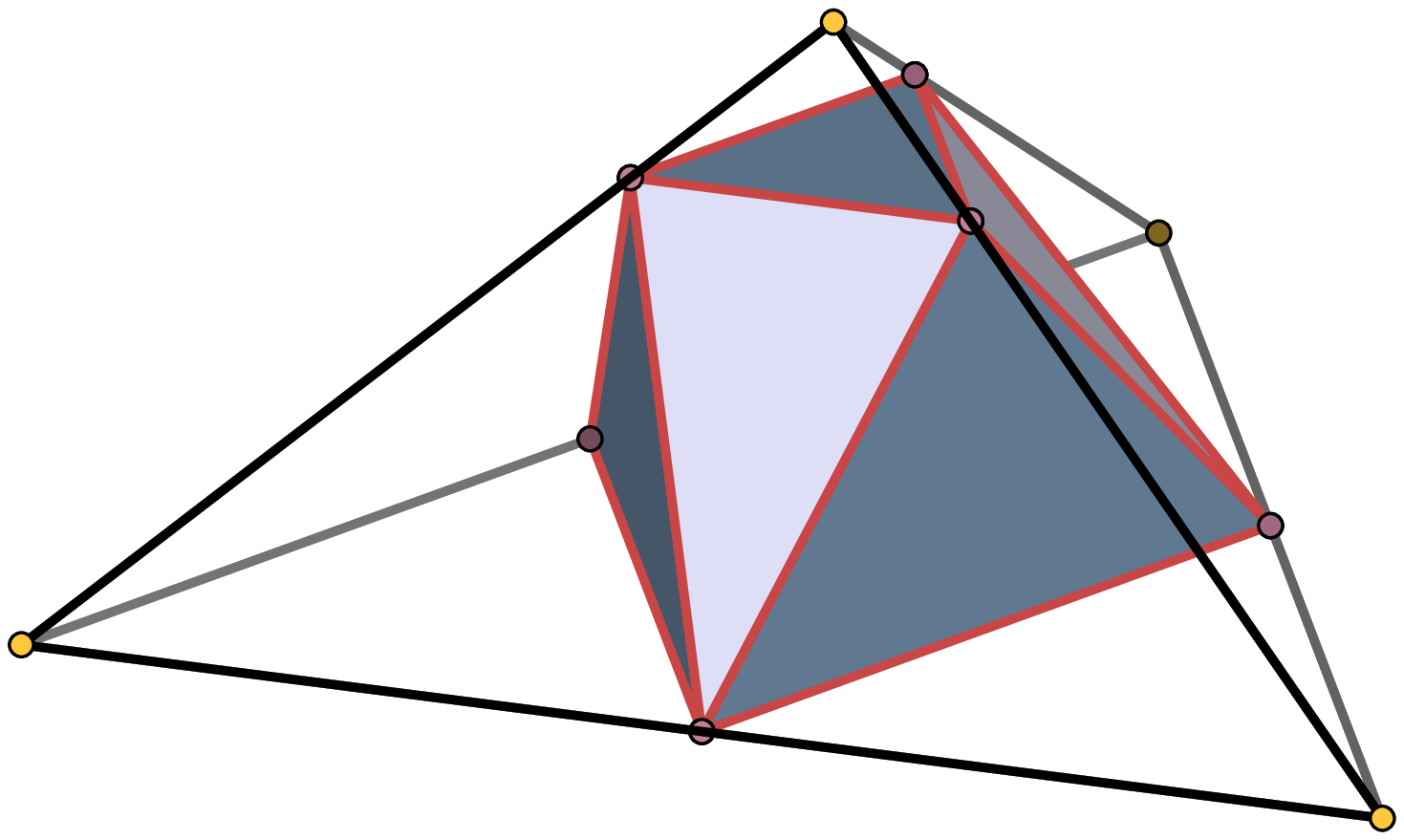} \hspace{0mm}
    \includegraphics[height=4cm,bb=36 111 469
    418,clip]{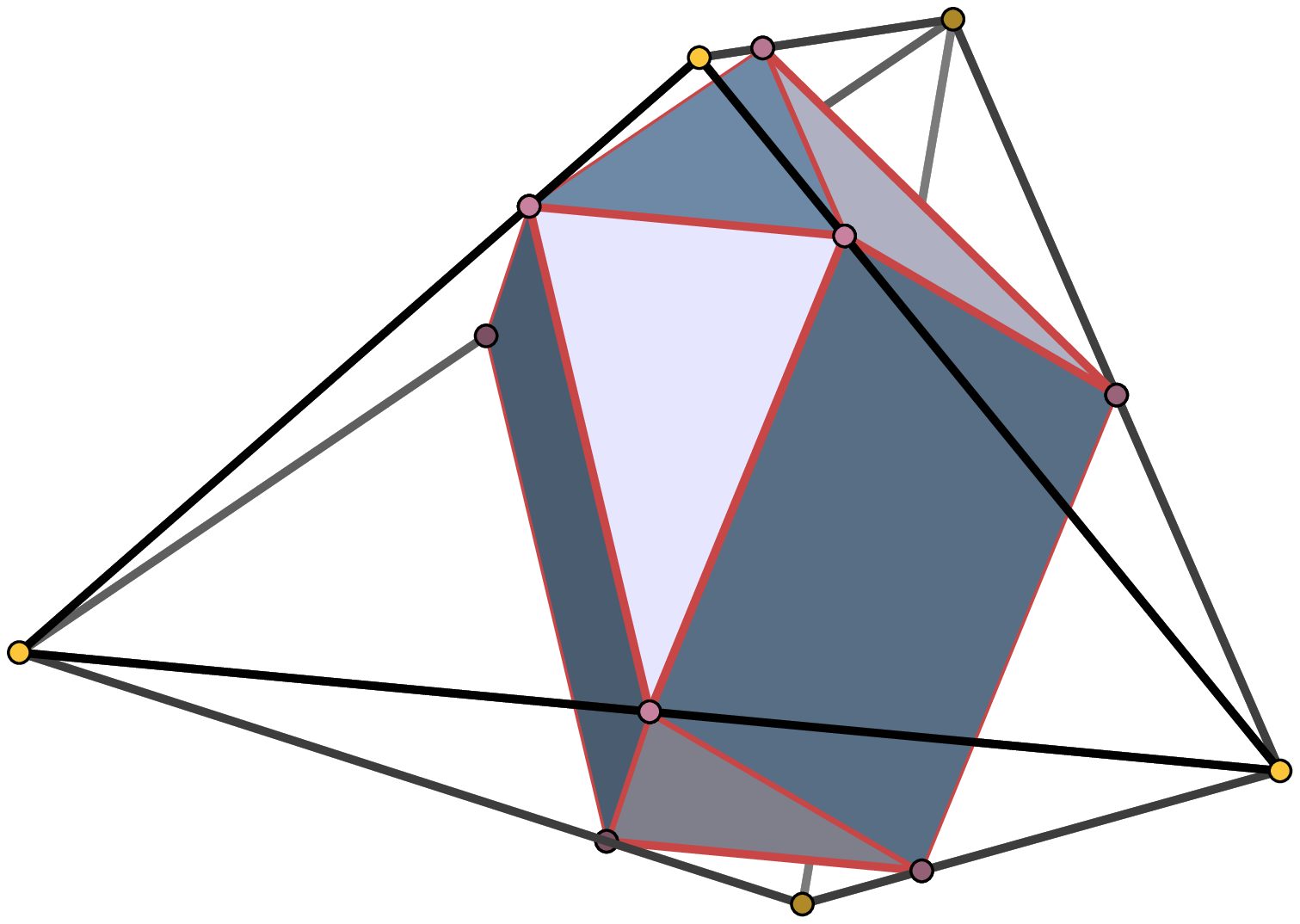} \hspace{0mm}
    \includegraphics[height=4cm,bb=57 117 480
    453,clip]{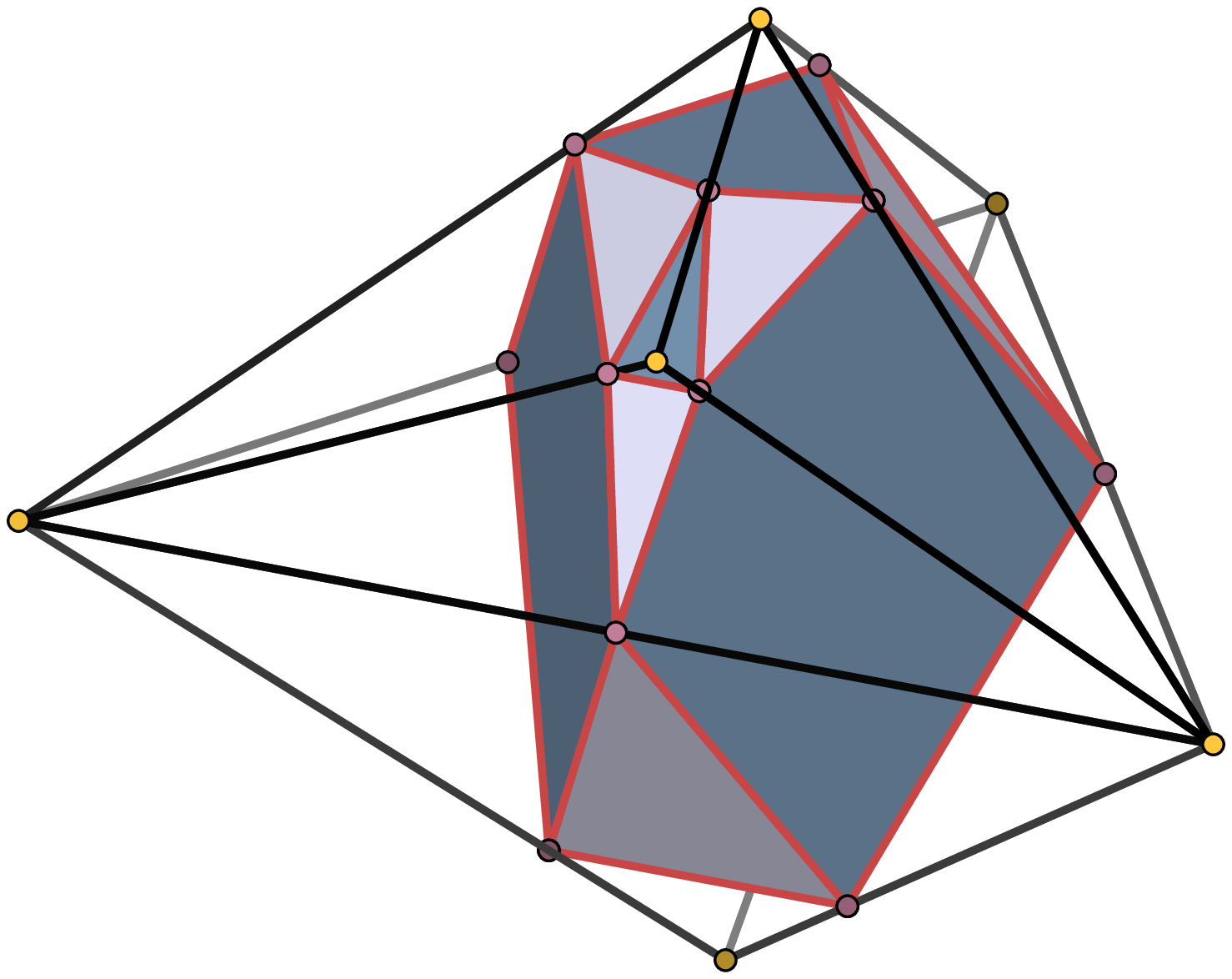} \hspace{-12mm}}
    \caption{Truncations of $P^3_0$, $P^3_1$, and $P^3_2$}
\end{figure}

In the construction for Theorem~\ref{thm:stacktrunc}, we can choose
the vertices of $P^d_0$ to be rational and take cutting hyperplanes
with rational normal vector. Further we can in all steps choose the
vertices we add to have rational coordinates. Thus also the
intersection points of the cutting hyperplanes with the edges will be
rational and we obtain the following corollary.

\begin{corollary}\label{cor:D1P}
  There are infinitely many combinatorially distinct rational
  $2$-simplicial $2$-simple $4$-polytopes $D_1(P^4_n)$, $n\ge0$, with
  the essential components $(f_0,f_1,f_2,f_3;f_{03})$ of the flag
  vectors given by
\[
f(D_1(P^4_n))\ \ =\ \ (10+4n,30+18n,30+18n,10+4n;50+26n).
\]\vskip-11mm
\qed
\end{corollary}

According to Eppstein, Kuperberg, and Ziegler \cite[Prop.~8]{Z80}, the
simplicial $4$-polytope $P^4_n$ has an edge-tangent realization if and
only if $n\le1$; this demonstrates that our vertex cutting approach is
indeed much more flexible than the approach via edge-tangent
realizations in \cite{Z80}.  The examples that we reproduce here are
the hypersimplex for $n=0$, and the glued hypersimplex of Braden
\cite{Bra-97}, which we get as $D_1(P^4_1)$. 

A similar infinite sequence of rational $2$-simple, $2$-simplicial
$4$-polytopes may be obtained from a stack of $n\ge1$
cross-polytopes. Using appropriate coordinates we obtain a realization
that is rational and has the symmetries of a regular $3$-simplex.
Thus we obtain a simplicial polytope $C_n^4$, with flag vector
$f(C^4_n)=(4+4n, 6+18n, 4+28n,2+14n; 8+56n)$.  Using the symmetry of
the polytope and induction on~$n$, it is easy to verify that $C^4_n$
in a suitable realization is indeed truncatable.

\begin{corollary}\label{cor:D1C}
  There is an infinite sequence of rational $2$-simplicial $2$-simple
  $4$-polytopes $D_1(C^4_n)$, $n\ge1$, with flag vectors
\[
f(D_1(C^4_n))\ \ =\ \ (6+18n,12+84n,12+84n,6+18n;24+120n).
\]
\end{corollary}

\begin{figure}[htb]
    \begin{center}
      \includegraphics[height=5cm]{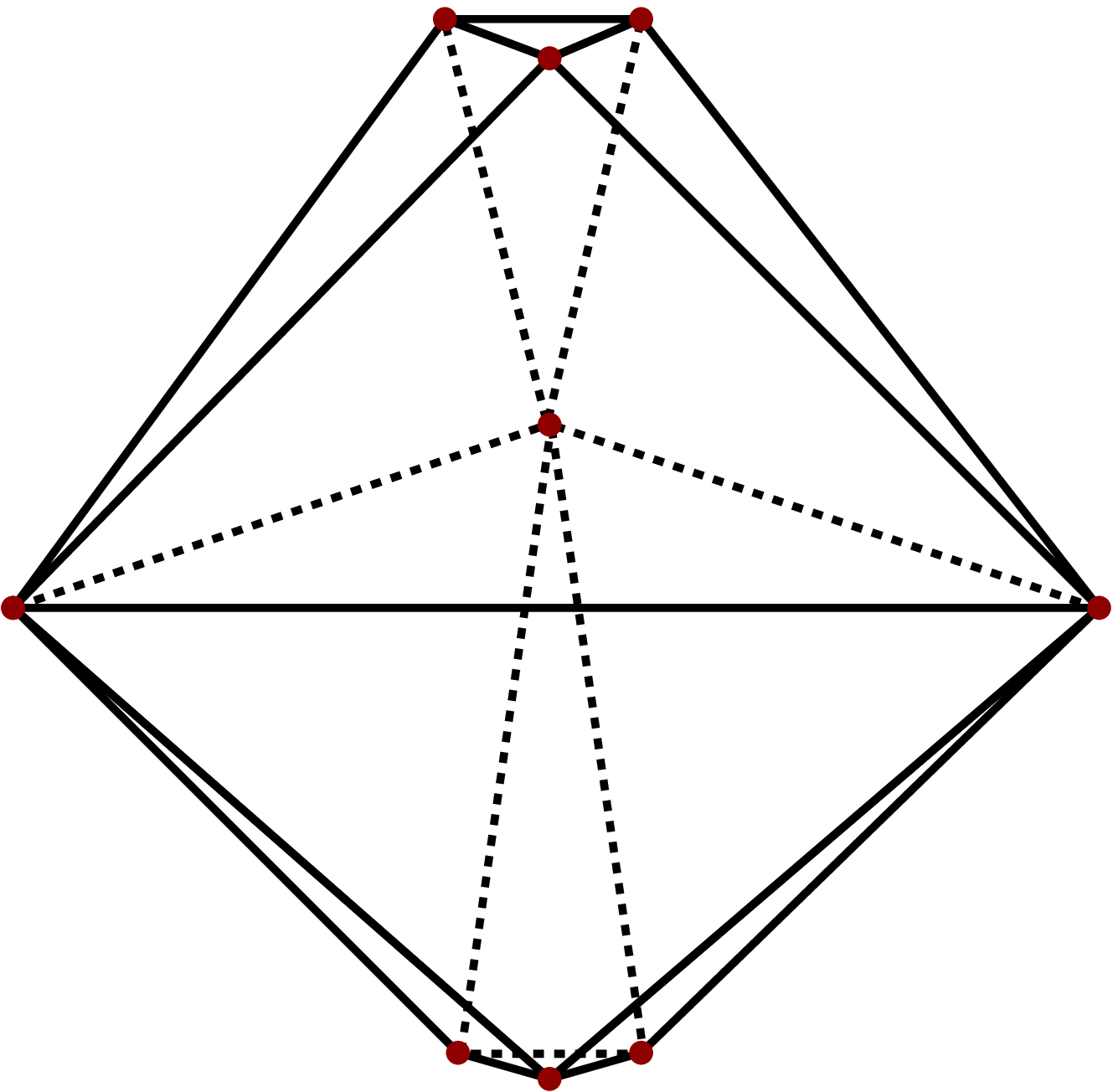}
      \hspace{12mm}
      \includegraphics[height=5cm]{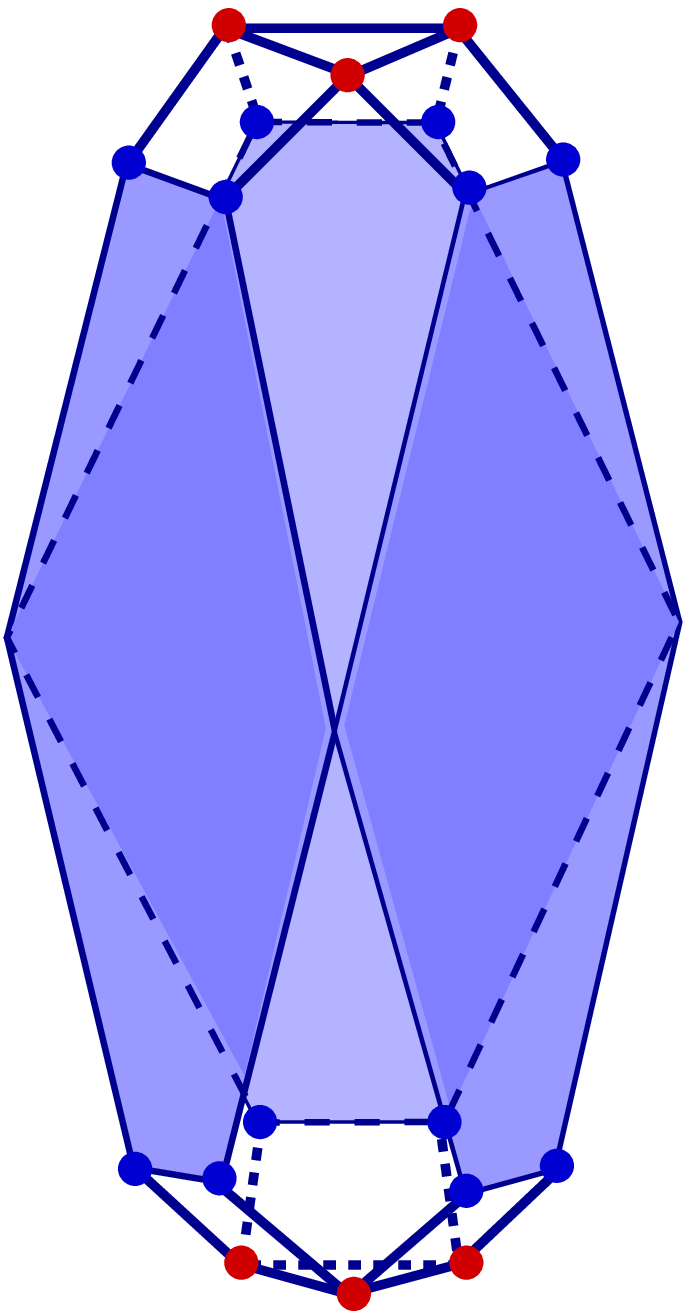}
      \hspace{12mm} \includegraphics[height=5cm,bb=200 90 435
      420,clip]{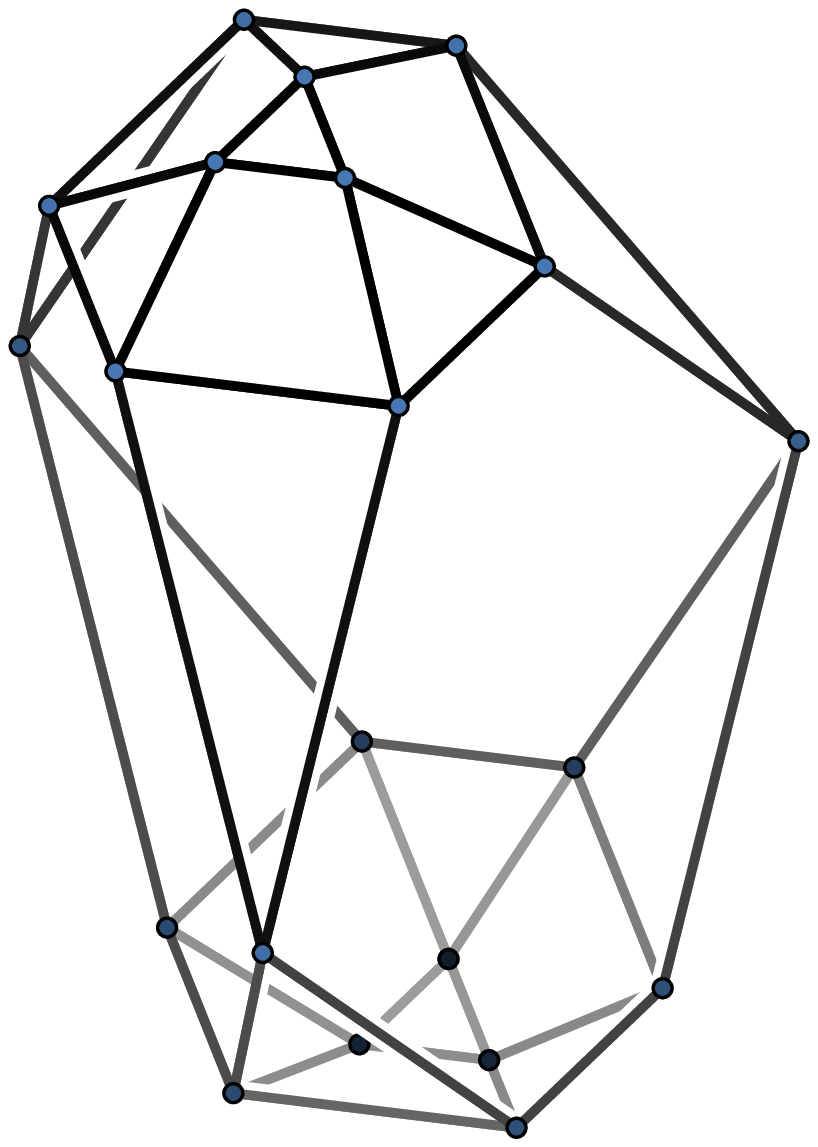}
    \end{center}
  \caption{Construction of $D_1(C^3_2)$: The left figure shows
    $C_2^3$; the center figure is obtained when three vertices have
    been truncated; the right figure displays $D_1(C_2^3)$}
  \label{C2}
\end{figure}

We leave the formal proof (and the explicit construction) to the
readers, and instead present graphics for the construction of
$D_1(C^3_2)$, in Figure~\ref{C2}.

\subsection{Realizations via hyperbolic geometry}

Now we extend the $E$-construction of \cite{Z80} to higher dimensional
polytopes. In dimension $d=4$ and for $t=1$ our construction will
coincide with this $E$-construction. The main goal of this section is
to prove the following result.

\begin{theorem}\label{thm:HigherDimensions}
  For every $d\ge3$ there are infinitely many combinatorially distinct
  $(d-2)$-simplicial $2$-simple $d$-polytopes.
\end{theorem}

The construction will roughly be as follows. First we prove that for
any $d$-polytope that has its $t$-faces ($0\le t\le d-1$) tangent to
the unit sphere $\S^{d-1}$ we can obtain a geometric realization of
$E_t(P)$ by taking the convex hull of the vertices of $P$ and its
dual. This is interesting mainly for the case $t=d-3$ and simplicial
$d$-polytopes, when Corollary \ref{cor:2s2s} tells us that
$E_{d-3}(P)$ will be $(d-2)$-simplicial and $2$-simple. In a next step
we will prove that there are in fact infinitely many combinatorially
different simplicial $d$-polytopes in any dimension $d\ge4$ that have
their $(d-3)$-faces tangent to $\S^{d-1}$.

The theorem is trivial in dimension $d=3$, as there are infinitely
many simple $3$-polytopes. All but the last step of our construction
will also work in dimension $3$, so from now on we will focus on
dimensions $d\ge4$, but all illustrations refer to dimension $3$. 

\begin{theorem}\label{thm:FaceTangentRealization}
  Let $P$ be a $d$-polytope and $1\le t\le d-2$. If $P$ has its
  $t$-faces tangent to the unit sphere $\S^{d-1}$ and $P^*$ is the
  polar of $P$, then $Q:=\conv(P\cup P^*)$
  is a polytopal realization of~$E_t(P)$.
\end{theorem}

\begin{proof}
  If a $t$-face $F$ of $P$ is tangent to $\S^{d-1}$ in a point $x\in
  F$, then all vertices of $F$ lie in the tangent space $T_x\S^{d-1}$
  of $\S^{d-1}$. By the definition of the polar polytope $P^*$ of $P$,
  the face $F^\circ$ of $P^*$ dual to $F$ is also contained in
  $T_x\S^{d-1}$ and is orthogonal to $F$ in $T_x\S^{d-1}$. Thus the
  convex hull $B(F)$ of $F$ and $F^\circ$ has dimension $d-1$ and is
  tangent to $\S^{d-1}$ in $x$ (see Figure~\ref{tangent_sphere}).
  \begin{figure}[htb]
    \begin{center}
      \psfrag{F}{$F$}
      \psfrag{S}{$\S^{d-1}$} 
      \psfrag{FD}{$F^\circ$}
      \psfrag{TSx}{\small $T_x\S^{d-1}$}
      \psfrag{EF}{$B(F)$}
      \psfrag{x}{$x$}
      \includegraphics[width=6cm]{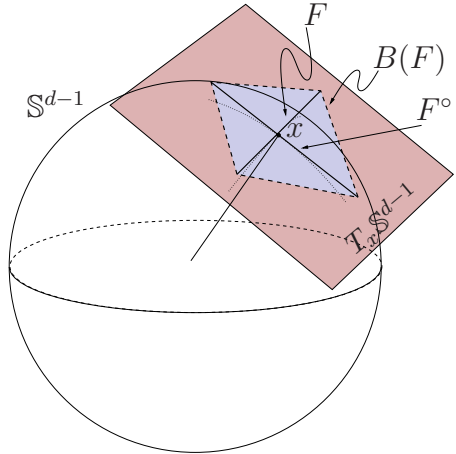}
      \caption{A $t$-face and its polar}
      \label{tangent_sphere}
    \end{center}
  \end{figure}
  From this we see that $Q=\conv(P\cup P^*)$ has the orthogonal sums
  $B(F):=\conv(F\cup F^\circ)$ as facets. 

  There is a facet $B(F)$ of $Q$ corresponding to each $t$-face $F$
  of~$P$. The vertices of $B(F)$ are given on the one hand by the
  vertices of $F$, on the other hand by the vertices of $F^\circ$,
  which are dual to the facets of $P$ that contain~$F$.
  
  Finally, one has to check that all facets of~$Q$ are sums of
  the type~$B(F)$. For this we show that all facets sharing a ridge
  with a facet $B(F)$ for a $t$-face $F$ of $P$ are of the type
  $B(F')$ for some $t$-face $F'$ of $P$. The facets of $B(F)$ are the
  convex hulls of a facet of $F$ and a facet of $F^\circ$, that is, of
  a $(t-1)$-face $R_F\subset F$ of $P$ and of a $(d-t-2)$-face of
  $P^*$ that is dual to a $(t+1)$-face $R^F\supset F$ of $P$. As the
  face lattice of $P$ is Eulerian, there is precisely one other
  $t$-face $F'$ of $P$ in the interval $[R_F,R^F]$ of length $2$.
  $B(F')$ also has $\conv(R_F\cup R^F)$ as a ridge, by construction.
  
  Thus $Q$ has a facet $B(F)$ for every $t$-face $F$ of $P$ (and no
  other facets) and it has the same vertex-facet-incidences as
  $E_t(P)$.
\end{proof}

In view of Theorem \ref{thm:FaceTangentRealization}, to prove Theorem
\ref{thm:HigherDimensions} we need to construct infinitely many
simplicial polytopes whose $(d-3)$-faces are tangent to the unit
sphere $\S^{d-1}$. This is an extension of arguments that were given
in \cite{Z80}, so we will only sketch the proof and point out the
differences.

We will view the interior $\D^d$ of the unit sphere $\S^{d-1}$ as the
the Klein model of hyperbolic space $\H^d$, and the unit sphere
$\S^{d-1}$ itself as the boundary at infinity $H^d_\infty$.
(See e.\,g.\ \cite[p.~75]{Iversen92} or
\cite{Ratcliffe1994}.)  The
advantage of this model of hyperbolic space in our situation is that
hyperbolic hyperplanes are the intersections of the ball with
Euclidean affine planes.

Now, for simplicity, we define a \emph{$T$-polytope} to be a
simplicial $d$-polytope that has its $(d-3)$-faces tangent to the unit
sphere $\S^{d-1}$.  The facets and ridges of a $T$-polytope $P$
properly intersect $\H^d$, the $(d-3)$-faces are tangent to
$\S^{d-1}\,\widehat=\,\H^d_\infty$ and all lower dimensional faces lie
outside the sphere. Thus the intersection $P^{hyp}:=\H^d\cap P$ is a
convex unbounded hyperbolic polyhedron in $\H^d$.

We repeat here a caveat from \cite{Z80}: A hyperideal hyperbolic
object -- even a convex polytope -- can be positioned in such a way
that it is unbounded as an Euclidean object (cf.
\cite[p.~508]{Schulte87}). However, we have the following lemma which
can be generalized to our situation (see also \cite{Springborn04}):
\begin{lemma*}[{\cite[Lemma 6]{Z80}}]
  For any edge-tangent convex polytope $Q$ in $\R^d$ whose points of
  tangency do not lie in a hyperplane there is a hyperbolic isometry
  $h$ whose extension to $\R^d$ maps $Q$ into a bounded position.\qed
\end{lemma*}

We want to compute the dihedral angle between any two adjacent facets
of a $T$-polytope $P$. As ridges and facets do at least partially lie
inside $\H^d$, this angle is well defined as a hyperbolic dihedral
angle of the hyperbolic polyhedron $P^{hyp}$ in $\H^d$ and is strictly
between $0$ and $\pi$.

The regular simplex and the regular cross polytope can be scaled so
that their $(d-3)$-faces are tangent to the unit sphere. Let
$\Delta_d^{hyp}$ resp.\ $C_d^{hyp}$ denote such realizations.

\begin{lemma}
  The dihedral angle between any two adjacent facets of
  $\Delta_d^{hyp}$ resp.\ $C_d^{hyp}$ is~${\pi}/{2}$
  resp.~${\pi}/{3}$.
\end{lemma}

\begin{proof}
  This is straightforward from the arguments given in \cite{Z80} when
  you keep in mind that their edges are in fact (codimension
  $3$)-faces for the purpose of our proof.
\end{proof}

\begin{figure}[htb]
  \begin{center}\qquad
    \subfigure[Two glued cross polytopes, with bipyramid
    facets]{\qquad \includegraphics[bb=153 104 490 407,
      clip,height=4cm]{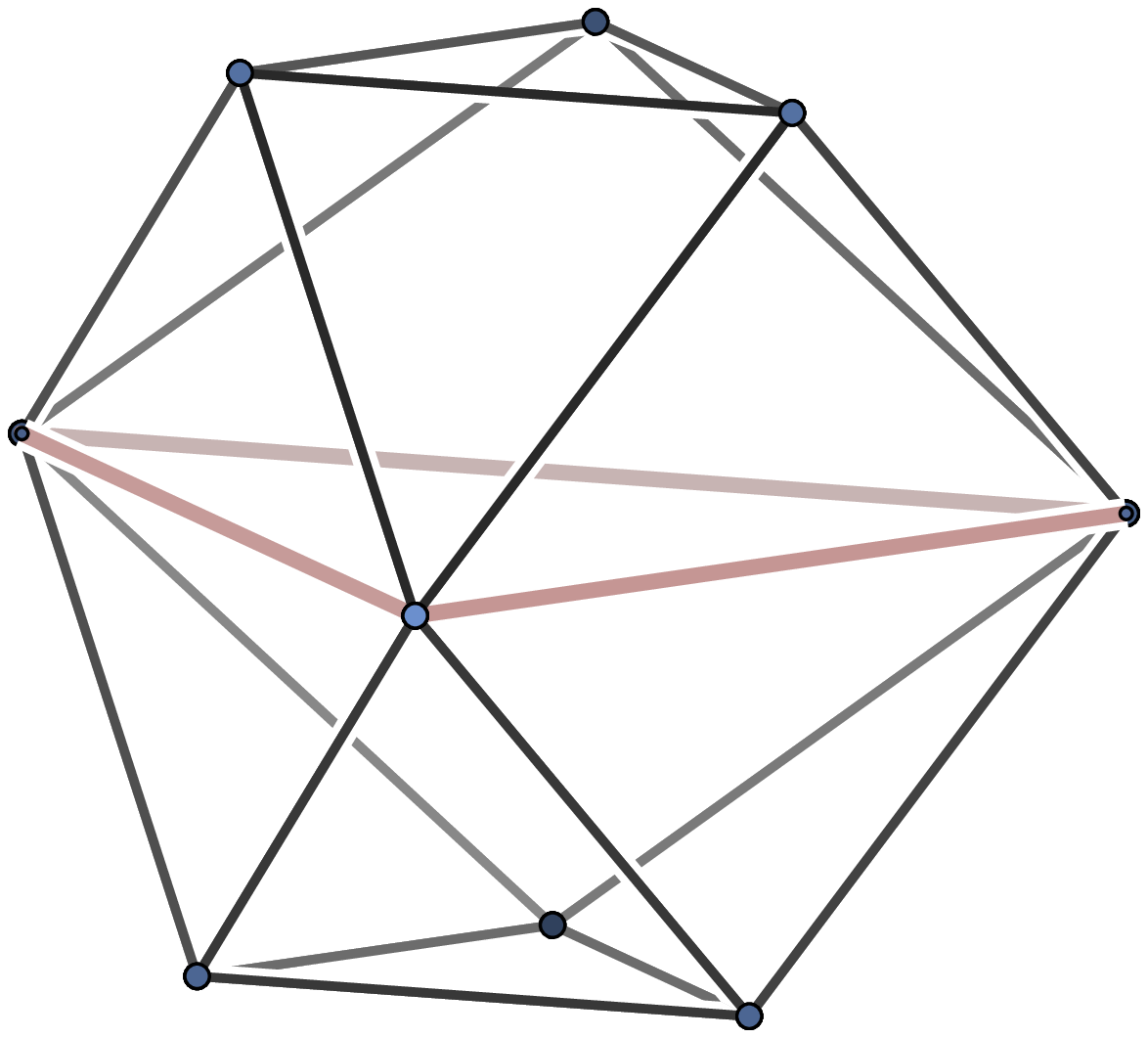}
      \label{fig:BadLink}\qquad}
    \subfigure[The link, with added simplices]{\qquad
      \psfrag{ls}{$\text{\footnotesize link of}
        \atop\mbox{\footnotesize a simplex}$}
      \psfrag{lc}{${\text{\footnotesize link of}
          \atop\mbox{\footnotesize a cross polytope}}$}
      \includegraphics[height=4cm]{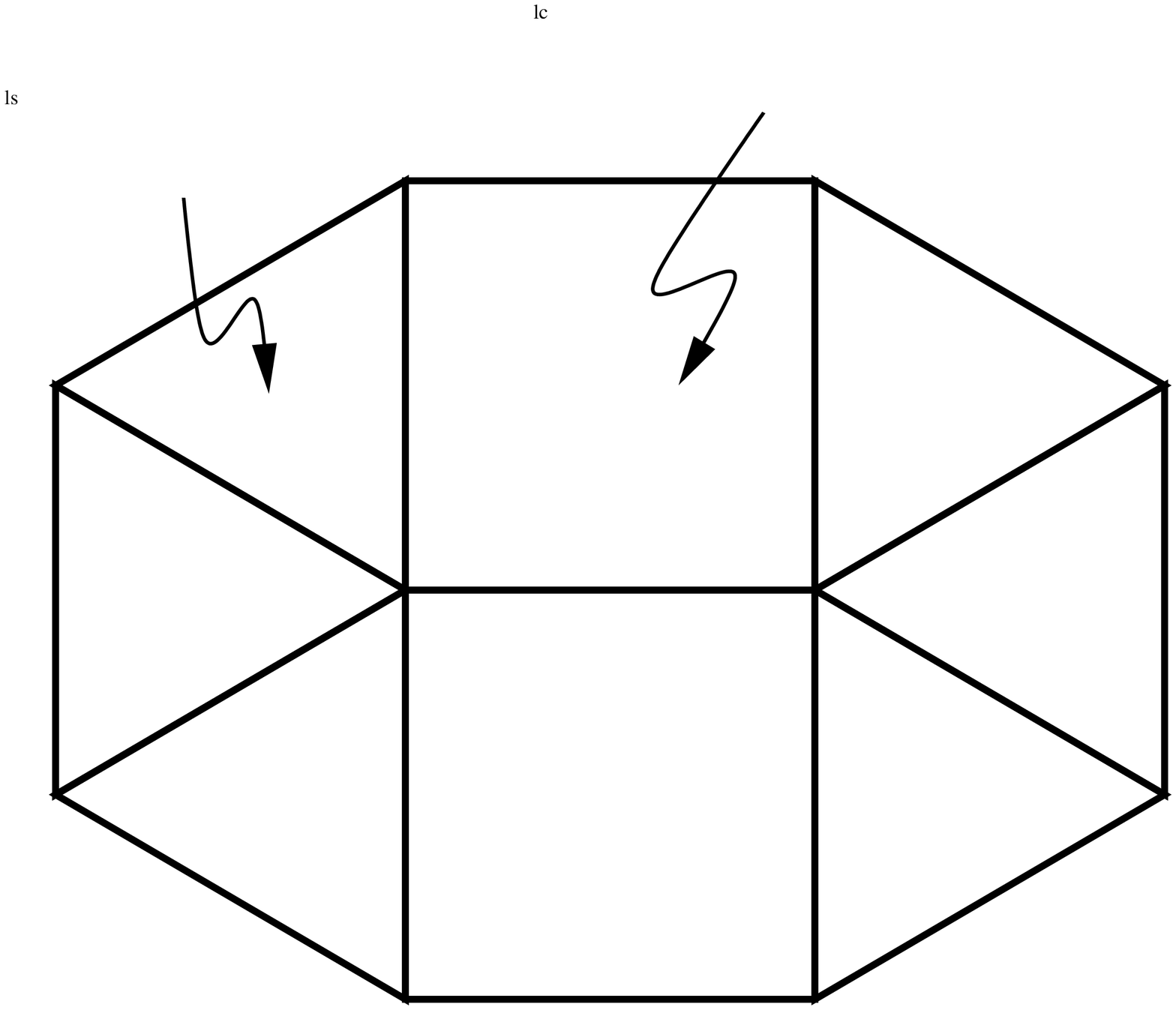}
      \label{link}\qquad\qquad}
   \subfigure[Three hyperideal tetrahedra have a bipyramidal facet (the top facet)]
      {\qquad \includegraphics[height=4cm]{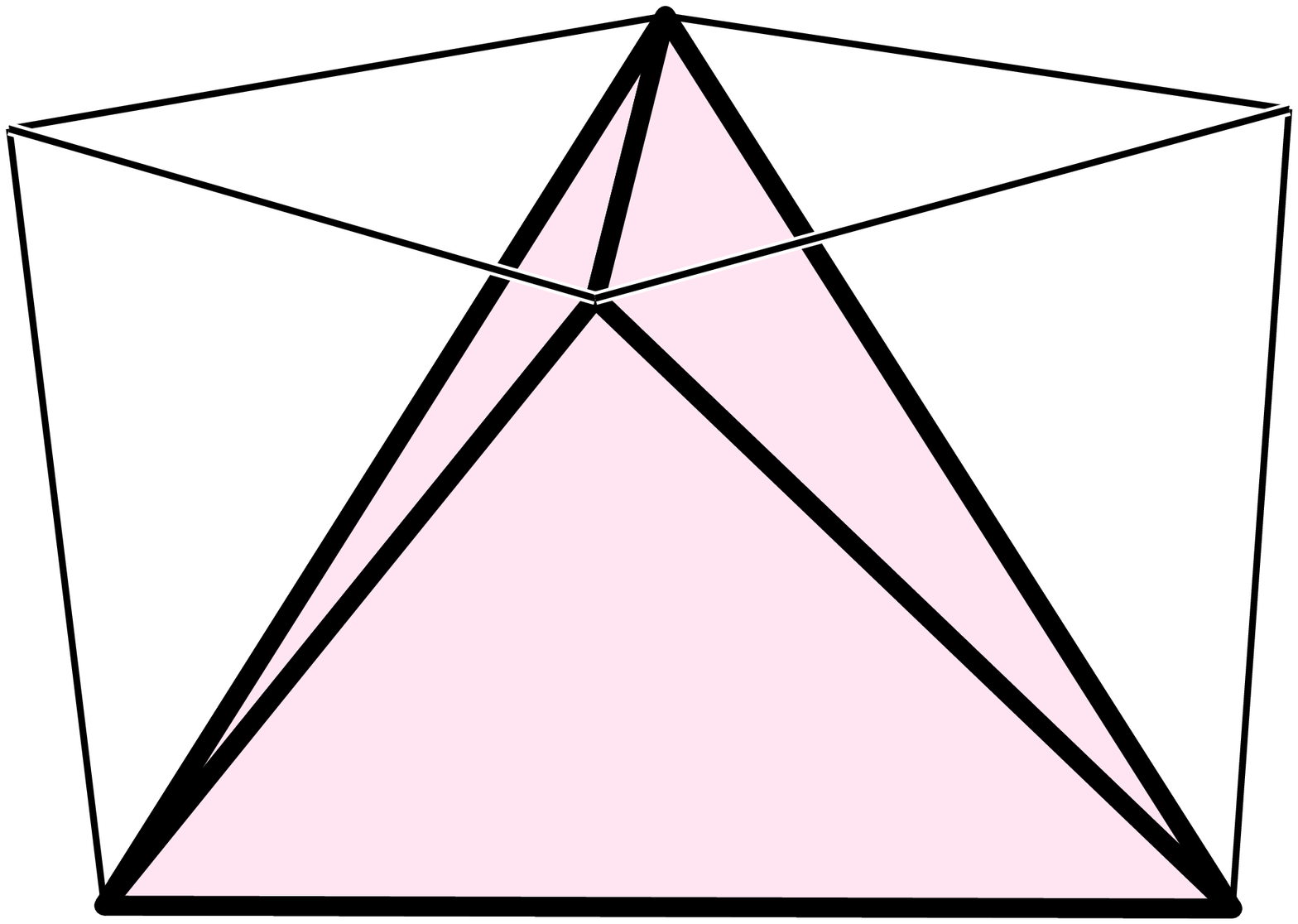} 
      \label{3tetra}\qquad\qquad}
    \caption{The gluing construction}
  \end{center}
\end{figure}

\begin{proof}[Proof of Theorem \ref{thm:HigherDimensions}]
  We will glue together certain $T$-polytopes along facets to get new
  $T$-polytopes. For this we have to position the $T$-polytopes in
  such a way that the two facets in question match, they both remain
  $T$-polytopes, and such that the dihedral angle between any two
  facets that become adjacent after the gluing is less than $\pi$
  (that is, the resulting compound should be convex).
  
  The isometry group of $\H^d$ is transitive and thus we can always
  map a facet of a regular simplicial $T$-polytope onto any facet of
  any other regular simplicial $T$-polytope. Thus gluing two such
  polytopes is always possible (when interpreted as hyperbolic
  objects).
  
  Now we glue together two cross polytopes along a facet $F$. As the
  dihedral angle of two adjacent facets of a cross polytope is
  $\frac{\pi}{2}$, the resulting polytope will have facets that are
  bipyramids (compare the shaded ridges in the stack of two glued
  cross polytopes drawn in Figure~\ref{fig:BadLink}), but we can
  repair this by gluing three simplices along any of these ridges in a
  way illustrated in Figure~\ref{link}. See Figure \ref{3tetra}
  for an example of three hyperideal simplices and the bipyramidal
  facet which we use for the gluing. The dihedral angle between two
  simplices is $\frac{2\pi}{3}$, and between a simplex and a cross
  polytope facet it is $\frac{5\pi}{6}$. Two such triples of glued-in
  simplices do not share a ridge if we are in dimension $d\ge4$.
  
  Iterating this construction by gluing $n$ cross polytopes
  ``end-to-end'' and gluing triples of simplices to all ridges having
  straight adjacent facets yields simplicial $T$-polytopes $Q_n^d$ with
  $f$-vectors
\[
f_j(Q_n^d)\ \ =\ \ 
      \begin{cases}
        2^dn-2(n-1)+d(3d-5)(n-1)              &\textrm{for } j=d-1,\\
        2^{d-1}dn-d(n-1)+\tfrac12 d^2(3d-5)(n-1)   &\textrm{for } j=d-2,\\
        2^{j+1}\binom{d}{j+1}n-\binom{d}{j+1}(n-1)
        +d\bigl\{3\binom{d}{j}-\binom{d-1}{j}\bigr\}(n-1) &\textrm{for
        } 0\le j\le d-2.
      \end{cases}
\]
for $n\ge1$. From these we then derive $(d-2)$-simplicial $2$-simple
$d$-polytopes $E_{d-3}(Q^d_n)$ whose $f$-vectors are given by equation
(\ref{eq:fvector}) as (recall $d\ge 4$)
\[
f_k(E_{d-3}(Q_n^d))\ \ =\ \ 
      \begin{cases}
        f_{d-3}(Q_n^d)
        &\textrm{for } k=d-1,\\[2mm]
        \binom{d-1}{2} f_{d-2}(Q^d_n) &\textrm{for } k=d-2,\\[2mm]
        \binom{d-1}{3} f_{d-2}(Q^d_n)
           +\binom{d}{3} f_{d-1}(Q^d_n) &\textrm{for } k=d-3,\\[2mm]
        \binom{d-1}{d-k} f_{d-2}(Q^d_n)+
           \binom{d}{d-k} f_{d-1}(Q^d_n)
            +f_{k}(Q^d_n) &\textrm{for } 1\le k\le
        d-4,\\[2mm]
        f_{d-1}(Q^d_n)+f_{0}(Q^d_n) &\textrm{for } k=0.
      \end{cases}
\]
For $d=4$ this formula indeed specializes to the $f$-vectors
$(54n-30,252n-156,252n-156,54n-30)$ of the $4$-dimensional examples
that were already constructed and computed in \cite{Z80}.
\end{proof}

The polytopes $E_{d-3}(Q_n^d)$ constructed in this proof have
nonrational coordinates, and there seems to be no easy way to remedy
this.


\subsection{Previously known examples}\label{subsec:KnownExamples}

Practically all examples of $(d-2)$-simplicial $2$-simple
$d$-polytopes ($d\ge4$) that appear in the literature may be seen as
special instances of the $E_t$-construction for spheres as presented
above, and realized by one of the two constructions of this section.
In particular, the Eppstein--Kuperberg--Ziegler \cite{Z80} examples
arise in the special case where $P$ is a simplicial $4$-polytope with
an edge-tangent realization $(t=1)$, or equivalently a simple
$4$-polytope with a ridge-tangent realization ($t=2$).  Most of the
other, earlier examples specialize even further to the case of regular
polytopes $P$, where the $t$-face tangency conditions can simply be
enforced by scaling --- but would typically yield irrational
coordinatizations, with no apparent degrees of freedom.

Until very recently, only finitely many $(d-2)$-simplicial $2$-simple
$d$-polytopes were known for every fixed $d\ge4$.  For $d=4$, in
addition to the simplex $\Delta_4$ this includes a few ``well-known''
examples:
\begin{compactitem}[$\bullet$]
\item Schl\"afli's self-dual $24$-cell, with flag vector
  $(24,96,96,24;144)$, arises as $E_2(4\textrm{-cube}) =
  E_1(4\textrm{-cross polytope})$.  In particular, a $24$-cell arose
  as $D_1(C^4_1)$ in Corollary~\ref{cor:D1C}.
\item The so-called hypersimplex is a $4$-polytope with
  vertex-transitive automorphism group, $5$ octahedra and $5$
  tetrahedra as facets, and flag vector $(10,30,30,10;50)$. It arises
  most easily as $D_1(\Delta_4)=D_1(P^4_0)$; see
  Corollary~\ref{cor:D1P}.  Its dual (with an automorphism group that
  is transitive on the $10$ bipyramid facets) thus is given by
  $E_1(\Delta_4)$.
\end{compactitem}
A much less obvious example is the glued hypersimplex of
Braden~\cite{Bra-97}, which arises as $D_1(P^4_1)$, by
Corollary~\ref{cor:D1P}.  There were claims by Perles and Shephard
(cf. \cite[pp.~82, 170]{Gr1-2}) for infinite families of
$2$-simplicial $2$-simple $4$-polytopes that turned out to be
premature.


\subsubsection{Gr\"unbaum's examples}

The $2$-simplicial $(d-2)$-simple $d$-polytopes ($d\ge4$) of
Gr\"unbaum \cite[pp.~65,66]{Gr1-2} are:
\begin{alignat*}{3}
  &K_k^d &&:=\bigl\{x\in\R^{d+1}\;\colon\; 0\le x_i\le 1\,,\,\sum_{i=0}^{d}
  x_i=k\bigr\} \qquad \mbox{for $1\le k\le d$}
  &&\qquad\text{\it ``hypersimplices''}\\
  &M^d &&:=\bigl\{x\in\R^d\;\colon\; |x_i|\le 1\,,\,\sum_{i=1}^d|x_i|\le
  d-2\bigr\}.&&
  \intertext{%
Since $K_k^d\cong K_{d+1-k}^d$, there are only $\lfloor\frac{d}{2}\rfloor$
combinatorial types of hypersimplices $K_k^d$ for each dimension~$d\ge4$.
    Gr\"unbaum also gives the following example of a
    $(d-3)$-simplicial $3$-simple $d$-polytope, $d\ge4$:}
  &N^d &&:=\bigl\{x\in\R^d\;\colon\;\sum\nolimits_{i=1}^d\eps_ix_i\le
  d-2\,,\ \eps_i=\pm1\,,\,\#\{\eps_i=1\}\mbox{
    odd}\bigr\}&&\qquad\text{\it ``half-cubes.''}
\end{alignat*}

$\Delta_d$ clearly has a geometric realization in which all $t$-faces
are tangent to the unit sphere for any $0\le t\le d-1$. Also $N_d$ in
the given realization has its $(d-2)$-faces tangent to a sphere. The
above theorem shows that we can apply the $E_t$-construction and an
easy check proves the following combinatorial equivalences:
\begin{alignat*}{2}
  (K_k^d)^*&\ \cong\ E_{k-1}(\Delta_d)&\textrm{that is,}&\quad
   K_k^d\ \cong \ D_{k-1}(\Delta_d) \qquad \mbox{for
    $1\le k\le d$},\quad\mbox{and}\\
  (M^d)^*&\ \cong\ E_{d-2}(N^d)&\quad \textrm{that is,}&\quad 
   M^d\ \cong \ D_{d-2}(N^d)\cong D_1(N^d)^*,
\end{alignat*} 
with the obvious generalization $D_k$ of vertex truncation that
preserves only one (relative interior) point from each $k$-face.  The
polytope $N^d$ itself cannot be a result of an $E_t$-construction as
for $d\ge 6$ the polytope $N^d$ and its dual are at least $3$-simple.


\subsubsection{The Gosset--Elte polytopes, and Wythoff's construction}

\begin{figure}[htb]
  \begin{center}
    \psfrag{r}{$\underbrace{\rule{2.5cm}{0cm}}_r$}
    \psfrag{s}{$\underbrace{\rule{2.5cm}{0cm}}_s$}
    \psfrag{t}{$\left\}\rule{0cm}{.6cm}\scriptstyle t\right.$}
    \includegraphics[width=6cm]{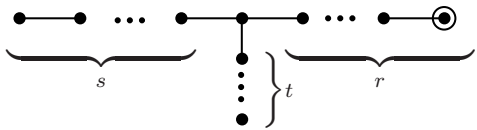}
    \caption{The Gosset--Elte polytopes}
    \label{coxeter}
  \end{center}
\end{figure}

Recently, Peter McMullen \cite{Gr1-2} noted that the Gosset--Elte
polytopes described in the classic book of Coxeter \cite{Cox} are
$(r+2)$-simplicial and $(s+t-1)$-simple. The Gosset--Elte polytopes
$r_{st}$ for $r,s,t\ge1$ arise from the Wythoff-construction in the
following way: Consider the group of reflections corresponding to the
diagram in Figure~\ref{coxeter}, where we have $r$ nodes on the right
end, $s$ nodes on the left end and $t$ nodes on the lower end. The
group of reflections is finite if and only if
$\frac{1}{r+1}+\frac{1}{s+1}+\frac{1}{t+1}>1$. This leaves us with
only three infinite series for $r,s$ and $t$ and a finite number of
other choices.

The three infinite series are $0_{d-k,k-1}$ for $1\le k\le d$,
$1_{d-3,1}$ and $(d-3)_{1,1}$. The first two coincide with
Gr\"unbaum's examples $K_k^d$ and $N^d$, while the last one gives the
$d$-dimensional cross polytope. Among the remaining finite number of
choices for $r,s$ and $t$, only $2_{2,1}$, $3_{2,1}$ and $4_{2,1}$ are
$(d-2)$-simplicial and $2$-simple. Out of these, only $0_{d-k-1,k}$
for $2\le k\le d-1$ is contained in our construction, as the
$1_{d-3,1}$ are $3$-simple and all others have simplicial facets,
which is not possible for a polytope resulting from the
$E_t$-construction for $0<t<d-1$.


\subsubsection{G\'evay's polytopes}

G\'evay constructs a number of interesting uniform polytopes as
``Kepler hypersolids'' that may be obtained as follows: Take $P$ to be
a regular polytope, scale it such that its $t$-faces are tangent to
the unit sphere, and then let $f_t(P)$ be the convex hull of $P$ and
its dual.  The polytopes thus obtained are dual to the corresponding
Wythoff/Gosset--Elte polytopes \cite[p.~128]{Gevay}.

As a specific instance, G\'evay \cite{Gevay,Gevay2} describes the
``dipyramidal $120$-cell'' $f_1H_4$, which appears in
Eppstein--Kuperberg--Ziegler \cite{Z80} as $E(120\textrm{-cell})$.


\subsubsection{Overview}

All examples except for those from \cite{Z80} are collected in Table
\ref{examples}. The first part of the table lists $4$-dimensional
examples, and in the second part examples in higher dimensions. In
dimension $d=3$ we have infinitely many $2$-simple (i.\,e.\ simple)
polytopes. The flag vector of a general $4$-polytope will be denoted
as $(f_0, f_1,f_2, f_3; f_{03})$. Any $2$-simplicial $2$-simple
$4$-polytope has flag vector of the special form $(f_0,f_1,f_1, f_0;
f_1+ 2f_0)$.

\begin{table}[htb]\label{examples}
{\footnotesize
\begin{tabular}{|l|l|l|l|l|l|}
\hline
 type             & flag vector & Gosset--Elte & Gr\"unbaum&G\'evay&here\\
\hline
simplex           &$(5,10,10,5;20)$ &$0_{3,0}=0_{0,3}$
                    &$K_4^4=K_1^4$  &$f_0A_4=f_3A_4$ &\\ 
hypersimplex      & $(10,30,30,10;50)$&
                    &&$f_1A_4=f_2A_4$ & $E_2(\textrm{simplex})$\\
hypersimplex$^*$  & $(10,30,30,10;50)$ & $0_{1,2}=0_{2,1}$
                    &$K_2^4=K_3^4$ &&$D_1(\textrm{simplex})$\\
24-cell           & $(24,96,96,24;144)$    &&$M_4$         
                  & $f_1B_4$ &$E_2(\textrm{4-cube})$\\
                  & $(720,3600,$ &&
                  & $f_1H_4$              &$E_2(\textrm{120-cell})$\\
                  &\;\;\;$3600,720;5040)$&&&&\\
glued             &&&&&\\
\; hypersimplex   & $(14,48,48,14;76)$ &&
                  &                    &$E_2(\textrm{prism over simplex})$\\
                  &&&&&\\
\hline            &&&&&\\
simplex           &&$0_{d-1,0}=0_{0,d-1}$& $K_1^d=K_d^d$ &&\\
hypersimplices    &&$0_{d-k-1,k}$& $K_k^d$&& $E_{k-1}^d(\textrm{simplex})^*$\\
$\; 2\le k\le d-1$&&&&&\\
                  &&& $M_d$   && $E_{d-2}(N_d)^*$\\
                  &&&&&$=E_2(d\textrm{-cube})^*$\\
                  &&$2_{2,1}$&&&\\
                  &&$3_{2,1}$&&&\\
                  &&$4_{2,1}$&&&\\
\hline
\end{tabular}
}
\caption{Examples of $E_t$-polytopes that were known before}
\end{table}
\medskip


\section{Some corollaries}\label{sec:fVectors}

While $2$-simple $2$-simplicial (``2s2s'' for the purpose of this
\newcommand\tsts{2s2s}%
section) $4$-polytopes seemed hard to construct until very recently
(with only finitely many examples known), we have now achieved quite
some flexibility. In particular, the vertex truncation method of
Section~\ref{subsec:vt} makes it easy to construct examples, say with
specific conditions on the resulting flag vectors.  To demonstrate
this, we here derive a sequence of corollaries.

The $f$-vector of a \tsts\ $4$-polytope is necessarily symmetric
($f_0=f_3$ and $f_1=f_2$), and it also determines the flag vector,
with $f_{03} = f_1+2f_2$.  Thus the $f$-vectors of \tsts\ 
$4$-polytopes have at most two independent parameters. Our first
corollary shows that we do indeed need two parameters.

\begin{corollary}\label{2s2sdiffedges}
  There are \tsts\ $4$-polytopes that have the same numbers of
  vertices (and facets), but that have different numbers of edges (and
  ridges).
\end{corollary}

\begin{proof}
  The $4$-dimensional regular cross polytope $C^*$ is clearly
  truncatable. Furthermore, we can stack $4$-simplices onto nine out
  of its 16 facets, so that the resulting simplicial $4$-polytope with
  $f$-vector $(17,60,86,43)$ is truncatable.  Thus we obtain a \tsts\ 
  $4$-polytope $D_1(C^*)$ with $f$-vector $(60,258,258,60)$.
  
  On the other hand, vertex truncation of a stack of three cross
  polytopes as in Corollary~\ref{cor:D1C} yields $D_1(C^4_3)$ with
  $f$-vector $(60,264,264,60)$.
\end{proof}

\begin{corollary}
  The number of combinatorially distinct \tsts\ $4$-polytopes with the
  $f$-vector $f(D_1(P^4_n))=(10+4n,30+18n,30+18n,10+4n;50+26n)$ grows
  exponentially with~$n$.
\end{corollary}

\begin{proof}
  There are exponentially many stacked $4$-polytopes with $n+5$
  vertices.  This already follows from the fact that there are
  exponentially many (unlabelled) trees of maximal degree~$5$ on $n+1$
  vertices.
  
  Thus we need to see that the combinatorial type of any stacked
  $4$-polytope $P^4_n$ can be reconstructed from its vertex truncation
  $D_1(P^4_n)$. (See Corollary \ref{cor:D1P}.)  The facets of
  $D_1(P^4_n)$ are on the one hand truncated simplices $F'$, which are
  octahedra, on the other hand the vertex figures $F_v$ of $P^4_n$,
  which are stacked. Furthermore, two of the octahedra $F'$ and $G'$
  are adjacent if and only if the corresponding facets $F$ and $G$ of
  $P^4_n$ are adjacent. Thus we get from $D_1(P^4_n)$ the dual graph
  of $P^4_n$, which determines the combinatorial type of $P^4_n$ by
  the reconstruction theorem of Blind and Mani
  \cite[Sect.~3.4]{Ziegler95}.
\end{proof}

One can tell from the flag vector whether a polytope is
$2$-simplicial, since this amounts to the condition $f_{02}=3f_2$, and
similarly for $2$-simplicity.  Our next corollary shows that there is
no similar criterion to derive this information from the $f$-vector.

\begin{corollary}\label{nonunique} 
  A \tsts\ and a not-\tsts\ $4$-polytope can have the same $f$-vector.
\end{corollary}

\begin{proof}
  Using a hyperbolic glueing construction for a stack of $n$
  600-cells, Eppstein, Kuperberg and Ziegler \cite[Sect.~3.3]{Z80}
  produced simplicial edge-tangent $4$-polytopes $Q_n$ with
  $f$-vectors $(106n+14,666n+54,666n+54,106n+14)$ and thus \tsts\ 
  $4$-polytopes $E_1(Q_n)$ with $f$ vector
  $(54+666n,240+3360n,240+3360n,54+666n)$.  We take this for $n=13$,
  that is, we get $f(E_1(Q_{13})) = (8712,43920,43920,8712)$ for a
  polytope that has lots of facets that are bipyramids over pentagons,
  and lots of ``regular'' vertices that are contained in exactly $12$
  such bipyramids, with dodecahedral vertex figure.  Now truncating
  $80$ such ``regular vertices'' and stacking pyramids over the
  resulting dodecahedral facets as well as onto $80$ other bipyramidal
  facets that are not involved in this results in a not-\tsts\ 
  $4$-polytope $E_1(Q_{13})'$ with $f(E_1(Q_{13})') = (
  10392,48280,48480,10392)$. 
  
  The \tsts\ $4$-polytope $D_1(C^4_{577})$ of Corollary~\ref{cor:D1C}
  has exactly the same $f$-vector.
\end{proof}

In the investigations of Eppstein, Kuperberg, and Ziegler~\cite{Z80},
\tsts\ $4$-polytopes were of interest as they provided examples of
$4$-polytopes for which the ``fatness'' parameter $\Phi(P) :=
\frac{f_1+f_2}{f_0+f_3}$ is particularly large. We now produce two
$4$-polytopes on the same number of vertices, one of them not-\tsts,
where the \tsts\ example is less fat.

\begin{corollary}\label{lessfat} 
  There is a \tsts\ $4$-polytope with the same number of vertices and
  facets as a $4$-polytope that is not \tsts\, but has fewer edges and
  ridges.
\end{corollary}

\begin{proof}
  The ``dipyramidal $720$-cell''
  $E:=E_2(120$-cell$)=D_1(600$-cell$)^*$ as discussed above has
  $f$-vector $(720,3600,3600,720)$. Now we perform some operations
  that destroy $2$-simplicity and $2$-simpliciality: We truncate two
  vertices with dodecahedral vertex figure, and stack pyramids on the
  resulting dodecahedral facets, and we also stack pyramids onto two
  bipyramidal facets.  Thus we get a new polytope $E'$ with $f$-vector
  $(762,3714,3714,762)$.
 
  On the other hand, vertex truncation applied to a stack of 42 cross
  polytopes yields the \tsts\ $4$-polytope $D_1(C^4_{42})$ with
  $f$-vector $(762,3540,3540,\allowbreak 762)$.
\end{proof}


\section{Open Problems}

\begin{compactenum}[1.]
\item Prove for some polytope $P$ and some~$t$ that $E_t(P)$ does not
  have a polytopal realization.
\item Give an example of some $E_t(P)$ that does not have a rational
  realization.  $E_2(\textrm{120-cell})$ is a candidate for this.
\item Modify the $E_t$-construction in order to also produce
  $(d-3)$-simplicial $3$-simple polytopes, such as the half-cubes
  $N^d$.
\end{compactenum}


\begin{small}

\end{small}


\end{document}